\documentclass[11pt]{article}
\usepackage{inputenc}
\usepackage[colorlinks=true, bookmarksnumbered=true, bookmarksopen=true,
bookmarksopenlevel=3, pdfstartview=FitH, linkcolor=cyan, pdfmenubar=true,
pdftoolbar=true, bookmarks=true,citecolor=cyan, urlcolor=magenta,
filecolor=cyan,menucolor=black,plainpages=false,pdfpagelabels]{hyperref}
\usepackage[lmargin=2cm,tmargin=2cm,bmargin=2cm,rmargin=2cm]{geometry}
\usepackage{lineno,hyperref}
\usepackage{amsfonts}
\usepackage{amsmath,amssymb}
\usepackage{mathtools}
\usepackage{color}
\usepackage{amsmath,amssymb}
\usepackage{amsthm}
\usepackage{amsfonts}
\usepackage{tikz}
\usepackage{pstricks}
\usepackage{pstricks-add}
\usepackage{pst-all,pst-math,pst-plot}
\usepackage{textcomp}
\usepackage{graphicx,color}
\usepackage{graphicx}
\usepackage{fancyhdr,ifthen}

\modulolinenumbers[5]
\mathtoolsset{showonlyrefs}
\usetikzlibrary{matrix}
\newtheorem{theorem}{Theorem}[section]
\newtheorem{lemma}[theorem]{Lemma}
\newtheorem{proposition}[theorem]{Proposition}
\newtheorem{corollary}[theorem]{Corollary}
\newtheorem{definition}[theorem]{Definition}
\newtheorem{remark}[theorem]{Remark}

\def\fin { \vskip 0pt \hfill $\diamond$ \vskip 12pt}

\numberwithin{equation}{section}

\begin{document}

\title{On the Helmholtz decomposition in Morrey and block spaces}
\author{{Lucas C. F. Ferreira$^{1}$}{\thanks{{Corresponding author. }\newline
{E-mail addresses: lcff@ime.unicamp.br (L.C.F. Ferreira),
marcosgabrieldesantana@gmail.com (M.G. Santana).}}} \ \ \ and \ \ {Marcos G.
Santana$^{2}$} \\
{\small $^{1,2}$ State University of Campinas, IMECC-Department of
Mathematics,} \\
{\small {Rua S\'{e}rgio Buarque de Holanda, 651, CEP 13083-859, Campinas-SP,
Brazil.}}}
\date{}
\maketitle

\begin{abstract}
In this work, we obtain the Helmholtz decomposition for vector fields in
Morrey, Zorko, and block spaces over bounded or exterior $C^{1}$ domains.
Generally speaking, our proofs rely on a careful interplay of localization,
flattening, and duality arguments. To accomplish this, we need to extend
some classical tools in analysis and PDE theory to those spaces, including
Stein extensions, compact embeddings, Poincar\'{e} inequalities,
Bogovskii-type theorem, among other ingredients. Some of these findings may
be of independent interest and applied to the study of a number of PDEs.


\end{abstract}

\tableofcontents


\renewcommand{\abstractname}{Abstract}

\pagestyle{fancy} \fancyhf{} \renewcommand{\headrulewidth}{0pt}
\chead{\ifthenelse{\isodd{\value{page}}}{Helmholtz decomposition in Morrey-type spaces}{L. C. F. Ferreira and M. G.
Santana}} \rhead{\thepage}

\section{Introduction}

Let $\Omega $ be a non-empty open subset of the whole space $\mathbb{R}^{n}$%
, where $n\geq 3$, and consider a vector field $\mathbf{u}:\Omega
\rightarrow \mathbb{R}^{n}$. The Helmholtz decomposition problem for $%
\mathbf{u}$ consists, formally, in finding a scalar function $p:\Omega
\rightarrow \mathbb{R}$ and a divergence-free vector field $\mathbf{v}%
:\Omega \rightarrow \mathbb{R}^{n}$ such that
\begin{equation}
\mathbf{u}=\nabla p+\mathbf{v}.  \label{decomposition intro}
\end{equation}%
If $\Omega =\mathbb{R}^{3}$ and $\mathbf{u}$ is smooth and decays
sufficiently fast at infinity, then \eqref{decomposition intro} holds
uniquely with
\begin{equation*}
p=\nabla \cdot \mathbf{U},\quad \mathbf{v}=-\nabla \times (\nabla \times
\mathbf{U}),
\end{equation*}%
where $\mathbf{U}$ is the solution of $\Delta \mathbf{U}=\mathbf{u}$. This
fact has been well known since H. von Helmholtz's work in electromagnetism
\cite{helmholtz1theorie}. Essentially, it states that vector fields can be
described by their irrotational (curl-free) and solenoidal (divergence-free)
components.

For more general classes of domains and function spaces, the Helmholtz
decomposition has been extensively studied for nearly a century due to its
crucial applications in important areas of physics and engineering, such as
elasticity, electromagnetism, and fluid mechanics. It has become a
fundamental tool in these fields, driving to the development of projection
methods in both analytical and numerical contexts (see, e.g., \cite%
{galdi2011introduction, Teman1}). Regarding fluid mechanics, let us
consider, for illustration, the formal heuristic for the well-known
Navier-Stokes equations
\begin{equation}
\partial _{t}\mathbf{w}-\Delta \mathbf{w}+(\mathbf{w}\cdot \nabla )\mathbf{w}%
+\nabla \phi =\mathbf{f},\quad \nabla \cdot \mathbf{w}=0,
\label{Navier-Stokes eq intro}
\end{equation}%
for the velocity field $\mathbf{w}:\Omega \times \lbrack 0,\infty
)\rightarrow \mathbb{R}^{n}$ and the scalar pressure $\phi :\Omega \times
\lbrack 0,\infty )\rightarrow \mathbb{R}$ of a viscous incompressible fluid
within a region $\Omega \subset \mathbb{R}^{n}$. For simplicity, we are
omitting the supplementary conditions. Suppose one is able to prove that %
\eqref{decomposition intro} holds both algebraically and topologically in a
certain Banach space $\mathcal{F}$ of vector fields $\mathbf{u}:\Omega
\rightarrow \mathbb{R}^{n}$; that is, we have a topological direct sum
decomposition of $\mathcal{F}$ induced by \eqref{decomposition intro}. Then,
there is a continuous projection $\mathbf{P}:\mathcal{F}\rightarrow \mathcal{%
F}$, the so-called \emph{Helmholtz projection}, defined by $\mathbf{P}%
\mathbf{u}=\mathbf{u}-\nabla p=\mathbf{v}$. Applying $\mathbf{P}$ to (\ref%
{Navier-Stokes eq intro}), we \ arrive at%
\begin{equation}
\partial _{t}\mathbf{w}-\mathbf{P}\Delta \mathbf{w}+\mathbf{P}(\mathbf{w}%
\cdot \nabla )\mathbf{w}=\mathbf{P}\mathbf{f}.
\label{abstract Navier-Stokes intro}
\end{equation}%
One can then search for a solution $\mathbf{w}$ to
\eqref{Navier-Stokes eq
intro} by studying \eqref{abstract Navier-Stokes intro} in the
divergence-free space $\mathbf{P}(\mathcal{F})$. A very powerful approach to
this is to investigate the properties of the \emph{Stokes operator} $-%
\mathbf{P}\Delta $, such as the existence of its fractional powers and
semigroup family \cite{Abe-Giga2013,
Giga-ken,bae2003analyticity,FUJITA-MORIMOTO1970,
Frohlich2007,Giga1981,giga1985domains,giga1988stokes}. The pressure $\phi $
can be recovered by taking the divergence in \eqref{Navier-Stokes eq intro}.

If $\Omega =\mathbb{R}^{n}$, then $\mathbf{P}$ can be described in terms of
the Riesz transforms $\mathcal{R}_{j}:=\partial _{j}(-\Delta )^{-1/2}$ as an
$(n\times n)$ matrix with elements $\delta _{ij}-\mathcal{R}_{i}\mathcal{R}%
_{j}$. In this case, the Helmholtz decomposition can be obtained in function
spaces where the Riesz transforms are bounded (see, e.g., \cite%
{Hobus-Saal-2019, kato1992strong, miyakawa1996hardy}). An analogous argument
applies if $\Omega $ is the half-space $\mathbb{R}_{+}^{n}$, as an explicit
formula for $\mathbf{P}$ can be derived using Green (Neumann) functions
(cf., e.g., \cite{galdi2011introduction}). For vector fields in Lebesgue
spaces $L^{q}(\Omega )$, with $1<q<\infty $, Fujiwara and Morimoto proved
the Helmholtz decomposition for bounded domains with smooth boundaries in
\cite{fujiwara1977theorem}. They established the existence of the normal
trace for vector fields in $L^{q}(\Omega )$ whose divergence also belongs to
$L^{q}(\Omega )$. Then, they combined this result with general $L^{q}$%
-theory results for boundary value problems of elliptic differential
equations. In \cite{simader1992new}, Simader and Sohr extended the
decomposition to bounded or exterior domains with $C^{1}$ boundaries (see
also \cite{miyakawa1982nonstationary}). To achieve this, they proved a
variational inequality (see also Theorem \ref{variational inequality theorem}%
) and applied it to the (weak) Neumann problem $\Delta p=\nabla \cdot
\mathbf{u}$ with the boundary condition $\partial p/\partial \mathbf{n}=%
\mathbf{u}\cdot \mathbf{n}$. In \cite{frohlich2000helmholtz}, Fr\"{o}hlich
generalized this approach to weighted Lebesgue spaces with Muckenhoupt
weights. For smooth exterior domains $\Omega \subset \mathbb{R}^{n}$,
Borchers and Miyakawa \cite{Borchers-Miyakawa1995} obtained the
decomposition in Lorentz spaces $L^{p,r}(\Omega )$ for $1<p<\infty $ and $%
1\leq r\leq \infty $, using a construction based on interpolation arguments.
By means of an approach based on suitable potential-type estimates for the
equivalent weak Neumann problem and interpolation arguments, Fujiwara and
Yamazaki \cite{fujiwara2007helmholtz} established the decomposition for
homogeneous Sobolev spaces $\dot{H}_{p}^{s}(\Omega )$ and Besov spaces $\dot{%
B}_{p,r}^{s}(\Omega )$ over bounded or exterior $C^{2,1}$ domains, where $%
p,r\in \lbrack 1,\infty ]$ and $1/p-1<s<1/p$ (with $p\neq 1,\infty $ in the
case of $\dot{H}_{p}^{s}$), as well as for their nonhomogeneous
counterparts. For arbitrary domains $\Omega $ of uniform $C^{2}$-type,
Farwig, Kozono, and Sohr \cite{Farwig-Kozo-Sohr2005, Farwig-Kozo-Sohr2007}
employed localization and covering arguments to obtain the Helmholtz
decomposition in the spaces $L^{p}(\Omega )\cap L^{2}(\Omega )$ ($2\leq
p<\infty $) and $L^{p}(\Omega )+L^{2}(\Omega )$ ($1<p<2$). In a sequence of
four recent papers, Giga and Gu established the Helmholtz decomposition in
the framework of $BMO$ spaces. In \cite{Giga-Gu-2020}, they used extension
and restriction-type arguments to prove the decomposition for $\mathbb{R}%
_{+}^{n}$. Later, they extended this result to bounded $C^{3}$ domains by
employing a potential-theoretic approach and deriving estimates for the
normal trace to solve the associated weak Neumann problem \cite{Giga-Gu-2021}%
. A similar approach was applied to slightly perturbed $C^{3}$ half-spaces
\cite{Giga-Gu-2023}. Additionally, they established the decomposition for
general uniformly $C^{3}$ domains by combining tools such as the Bogovskii
operator, Morrey and Poincar\'{e} inequalities in Sobolev and H\"{o}lder
spaces (localized appropriately with suitable control over the dependence on
constants), and constructing specific cut-off functions based on the
geometric properties of the domain \cite{giga2024helmholtz}.

In this work, we prove the Helmholtz decomposition for Morrey, Zorko and
block spaces (see Definitions \ref{definition:MorreyAndZorkoSpaces} and \ref%
{definition:BlockSpaces}) with $\Omega $ being either a bounded or exterior
domain with $C^{1}$ boundary, $\mathbb{R}^{n}$ or $\mathbb{R}_{+}^{n}$.
These spaces constitute generalizations of Lebesgue spaces, often arising
naturally in harmonic analysis and the study of partial differential
equations (PDEs), and are closely connected with other classical function
spaces \cite{adams2015morrey,
adams2012morrey,rosenthal2015morrey,sawano2020morrey, zorko1986morrey}. As a
matter of fact, Morrey spaces encompass a broader class than $L^{p}$ and
Lorentz spaces (with the same scaling), enabling the analysis of a wider
variety of functions. This flexibility is particularly useful for examining
initial-boundary data and solutions that exhibit varying degrees of
integrability and singularity.

The validity of the Helmholtz decomposition naturally motivates deeper
research on the already existing studies of the Navier-Stokes equations in
Morrey spaces \cite{deAlm-Fer2013, federbush1993navier, Fer1,
giga1989navier, kato1992strong, taylor1992analysis}. Moreover, the auxiliary
results developed here can provide useful tools for the study of this and
other PDEs. It is noteworthy that in \cite{kato1992strong}, Kato proved the
boundedness in Morrey spaces over $\Omega =\mathbb{R}^{n}$ of the projection
operator $\mathbf{P}$ constructed through Riesz transforms. For the sake of
completeness, however, we keep the proof of the Helmholtz decomposition for
this case here. Additionally, while Kato does not characterize the range and
kernel of $\mathbf{P}$ in his work, it is straightforward to verify that the
decomposition obtained here coincides with the one induced by his projection
operator.

One of the challenges in considering Morrey and block spaces is their lack
of suitable interpolation properties (see \cite{Blasco1, Fer1}), which
hinders the direct extension of results previously established in $L^{p}$ as
done in the case of Lorentz spaces. This structural difficulty requires
tailored approaches to analyze functions and operators in those spaces.
Moreover, unlike the case of Lebesgue spaces in \cite{simader1992new} or
Muckenhoupt weighted-Lebesgue spaces in \cite{frohlich2000helmholtz}, Morrey
spaces are non-reflexive (see also Remark \ref{remark weak-star closure})
and non-separable. Note that compactly supported smooth functions are not
dense in Morrey spaces (see Remark \ref{remark L infinite not dense}). To
overcome these obstacles, we first solve the decomposition in a less
singular subspace and then utilize duality relations to extend the result,
as briefly explained below.

Our proof of the Helmholtz decomposition is structured as follows. First, by
density arguments, we establish the decomposition for the Zorko space $%
\mathring{M}_{q,\lambda }\left( \Omega \right) $, $1<q<\infty $, $0\leq
\lambda <n$, which is defined as the closure of $C_{0}^{\infty }(\Omega )$
in the Morrey space $M_{q,\lambda }\left( \Omega \right) $. We utilize the
fact that this decomposition is equivalent to the well-posedness of the weak
version of the Neumann problem
\begin{align*}
\Delta p& =\nabla \cdot \mathbf{u},\quad \text{in }\Omega , \\
\frac{\partial p}{\partial \mathbf{n}}& =\mathbf{u}\cdot \mathbf{n},\quad
\text{on }\partial \Omega .
\end{align*}%
Furthermore, we divide the proof into cases based on the domain. The
simplest case, $\Omega =\mathbb{R}^{n}$ or $\mathbb{R}_{+}^{n}$, is handled
with the aid of Calder\'{o}n-Zygmund operators. For a bounded or exterior $%
C^{1}$ domain, we employ the technique of localization and flattening used
in \cite{simader1992new}. Subsequently, we establish the decomposition for
block spaces $H_{q,\lambda }(\Omega )$ and Morrey spaces through duality
arguments, considering the relations

\begin{equation*}
H_{q,\lambda }(\Omega )=\mathring{M}_{q^{\prime },\lambda }(\Omega )^{\ast }%
\text{ and }M_{q,\lambda }\left( \Omega \right) =H_{q^{\prime },\lambda
}\left( \Omega \right) ^{\ast }.
\end{equation*}

In order to perform the above construction, we need to extend some tools in
analysis and PDE theory to the framework of Morrey, Zorko, and block spaces,
covering topics such as mollification, approximation by smooth functions,
integral operators, Stein extensions, compact embeddings (Rellich
Kondrachov-type theorem), Poincar\'{e}\ inequalities, Bogovskii-type
theorem, among others. For example, we obtain a Stein-type extension in
Zorko spaces, which appears to hold independent interest. Here, such an
extension serves as a sort of first link in the chain of several steps in
the construction of the Helmholtz decomposition, being used particularly to
prove a Poincar\'{e} inequality in Zorko spaces. Other ingredients that may
be of interest in themselves are the Rellich Kondrachov-type embedding and
the Bogovskii-type theorem in our framework.

In what follows we give the precise statement of the Helmholtz decomposition
and projection, as well as the functional settings in which they are
addressed. Afterwards, for the reader convenience, we conclude the
introduction by providing a detailed description of the organization of the
manuscript.

\subsection{Statement of the main result}

\label{Sec-Statement}As already stated above, the domain $\Omega $ is a
non-empty open subset of $\mathbb{R}^{n}$, where the dimension $n\geq 3.$
When necessary, we outline additional conditions on $\Omega $, including
connectedness, smoothness, boundedness, and others. For instance, in our
main result, we will assume that $\Omega $ is either a bounded, an exterior
domain with a $C^{1}$ boundary, $\mathbb{R}^{n}$, or $\mathbb{R}%
_{+}^{n}:=\{x\in \mathbb{R}^{n}:x_{n}>0\}$).

We remind that $\Omega $ is called an exterior domain if $\mathbb{R}%
^{n}\backslash \Omega $ is compact. Moreover, unless explicitly noted, we
always suppose that the indexes $q$ and $\lambda $ belong to ranges $%
1<q<\infty $ and $0\leq \lambda <n$, respectively. We denote by $q^{\prime }$
the conjugated exponent of $q$.

For $R>0$ and $x\in \mathbb{R}^{n}$, we consider $B(x,R):=\{y\in \mathbb{R}%
^{n}:|x-y|<R\}$, and $\Omega _{R}(x):=\Omega \cap B(x,R)$ and $\Omega
^{R}(x):=\Omega \backslash \overline{B(x,R)}$. We abbreviate by $%
B_{R}:=B(0,R)$, $\Omega _{R}:=\Omega _{R}(0),\Omega ^{R}:=\Omega ^{R}(0)$.

\begin{definition}[Morrey and Zorko Spaces]
\label{definition:MorreyAndZorkoSpaces} We denote the (scalar or
vector-valued) Morrey space by $M_{q,\lambda}(\Omega)$, defined as the space
of all functions $f \in L^q_{\mathrm{loc}}(\overline{\Omega})$ such that
\begin{equation}  \label{norm of morrey}
\|f\|_{M_{q,\lambda}(\Omega)} := \sup \left\{ R^{-\lambda/q}
\|f\|_{L^q(\Omega_R(x))} : R > 0, \, x \in \Omega \right\} < \infty.
\end{equation}
The functional in \eqref{norm of morrey} defines a norm on $%
M_{q,\lambda}(\Omega)$, making it a Banach space.

We also introduce $\mathring{M}_{q,\lambda}\left(\Omega\right)$, known as
the \emph{Zorko space} (see \cite{zorko1986morrey}), which is defined as the
closure of $C_0^\infty(\overline{\Omega})$ in $M_{q,\lambda}(\Omega)$.
\end{definition}

We observe that, for $\lambda>0$, $\mathring{M}_{q,\lambda}\left(\Omega%
\right)$ is a proper subspace of $M_{q,\lambda}\left(\Omega\right)$ (see
Remark \ref{remark L infinite not dense}), while $M_{q,0}(\Omega)=\mathring{M%
}_{q,0}(\Omega)=L^q(\Omega)$.

\begin{definition}[Block Spaces]
\label{definition:BlockSpaces} Let $q^{\prime }$ denote the conjugate
exponent of $q$. A \emph{$(q,\lambda )$-block} is a function $\rho \in L_{%
\mathrm{loc}}^{q}(\overline{\Omega })$ such that, for some $R>0$ and $x\in
\Omega $, we have $\mathrm{supp}(\rho )\subset \overline{\Omega _{R}(x)}$
and
\begin{equation*}
R^{\lambda /q^{\prime }}\Vert \rho \Vert _{L^{q}(\Omega _{R}(x))}\leq 1.
\end{equation*}%
The \emph{block space} $H_{q,\lambda }(\Omega )$ is defined as the space of
all functions $\phi $ of the form
\begin{equation*}
\phi =\sum_{k=0}^{\infty }c_{k}\phi _{k},
\end{equation*}%
where each $\phi _{k}$ is a $(q,\lambda )$-block and $\{c_{k}\}\in \ell ^{1}$%
. The space $H_{q,\lambda }(\Omega )$ is a Banach space with the norm
\begin{equation*}
\Vert \phi \Vert _{H_{q,\lambda }(\Omega )}:=\inf \left\{ \Vert
\{c_{k}\}\Vert _{\ell ^{1}}:\phi =\sum_{k=0}^{\infty }c_{k}\phi _{k}\text{
with }\phi _{k}\text{ a }(q,\lambda )\text{-block}\right\} .
\end{equation*}
\end{definition}

Before state our main result, we remark the following duality relation
between Morrey, Zorko and block spaces, whose proof can be found in \cite%
{almeida2017approximation}:
\begin{equation}
\mathring{M}_{q,\lambda }\left( \Omega \right) ^{\ast }=H_{q^{\prime
},\lambda }(\Omega )\quad \text{and}\quad H_{q^{\prime },\lambda }(\Omega
)^{\ast }=M_{q,\lambda }\left( \Omega \right) ,
\label{equation:DualityRelations}
\end{equation}%
where the isomorphisms hold with respect to the duality induced by
integration and with norm equivalence.

\begin{definition}[Solenoidal and irrotational vector fields]
For $X=\mathring{M}_{q,\lambda}\left(\Omega\right)$, $H_{q,\lambda}(\Omega)$
or $M_{q,\lambda}\left(\Omega\right)$, we define the \emph{irrotational part
of $X$} by
\begin{equation*}
GX:=\{\nabla p : p \in W^{1,1}_{\mathrm{loc}}(\overline{\Omega})\}\cap X.
\end{equation*}
Also, we define the \emph{solenoidal part of $X$} by
\begin{equation*}
SX:=\text{ closure of } C^\infty_{0,\sigma}(\Omega) :=\{\mathbf{w} \in
C^\infty_0(\Omega)^n : \nabla\cdot \mathbf{w}=0\} \text{ in }X,
\end{equation*}
where the closure is taken in the strong topology if $X=\mathring{M}%
_{q,\lambda}\left(\Omega\right)$, and in their respective weak-star topology
\footnote{%
For Lebesgue spaces, that is, for $\lambda=0$, the spaces $SM_{q,0}(\Omega)$
and $SH_{q,0}(\Omega)$ defined above defined coincide with the \emph{strong}
closure of $C^\infty_{0,\sigma}(\Omega)$ in $L^q(\Omega)$, recovering the
usual definition in literature. See also Remark \ref{remark weak-star
closure}.} induced by \eqref{equation:DualityRelations} if $%
X=H_{q,\lambda}(\Omega)$ or $M_{q,\lambda}\left(\Omega\right)$.
\end{definition}

We are now ready to state our main result.

\begin{theorem}[Helmholtz Decomposition]
\label{theorem:HelmholtzDecomposition} Let $n\geq 3$ and $\Omega $ be either
a bounded, exterior domain with a $C^{1}$ boundary, $\mathbb{R}^{n}$, or $%
\mathbb{R}_{+}^{n}$. Let $1<q<\infty $, $0\leq \lambda <n$, and let $X$ be
either $\mathring{M}_{q,\lambda }\left( \Omega \right) $, $H_{q,\lambda
}(\Omega )$, or $M_{q,\lambda }\left( \Omega \right) $. Then:

The Helmholtz decomposition holds for $X$, i.e.,
\begin{equation*}
X=GX\oplus SX,
\end{equation*}%
as an algebraic and topological sum. This means that for each vector field $%
\mathbf{u}\in X$, there are unique $\mathbf{v}\in GX$ and $\mathbf{w}\in SX$
such that $\mathbf{u}=\mathbf{v}+\mathbf{w}$, and there exists a constant $%
c>0$ independent of $\mathbf{u}$ such that
\begin{equation*}
\Vert \mathbf{v}\Vert _{X}+\Vert \mathbf{w}\Vert _{X}\leq c\Vert \mathbf{u}%
\Vert _{X}.
\end{equation*}

In particular, there is a bounded projection (Helmholtz projection) $\mathbf{%
P}_X: X \to X$ with kernel $GX$ and range $SX$. Moreover, the following
duality relations hold:

\begin{itemize}
\item $SM_{q,\lambda} = \left(GH_{q^{\prime},\lambda}\right)^\perp$, $%
GM_{q,\lambda} = \left(SH_{q^{\prime},\lambda}\right)^\perp$, and $\mathbf{P}%
_{M_{q,\lambda}} = \left(\mathbf{P}_{H_{q^{\prime},\lambda}}\right)^\ast;$

\item $SH_{q,\lambda} = \left(G\mathring{M}_{q^{\prime},\lambda}\right)^%
\perp $, $GH_{q,\lambda} = \left(S\mathring{M}_{q^{\prime},\lambda}\right)^%
\perp$, and $\mathbf{P}_{H_{q,\lambda}} = \left(\mathbf{P}_{\mathring{M}%
_{q^{\prime},\lambda}}\right)^\ast.$
\end{itemize}
\end{theorem}

\begin{remark}
We note in advance that, as a consequence of the proof of the forthcoming
Theorem \ref{variational inequality theorem}, the Helmholtz decomposition
for Morrey, Zorko, and block spaces also holds when $\Omega $ is a $C^{1}$%
-smooth and slightly perturbed half-space, that is, if
\begin{equation*}
\Omega :=\{x=(x^{\prime },x_{n})\in \mathbb{R}^{n}:x_{n}>\sigma (x^{\prime
})\},
\end{equation*}%
for some function $\sigma \in C_{0}^{1}(\mathbb{R}^{n-1})$ such that $\Vert
\sigma \Vert _{L^{\infty }(\mathbb{R}^{n-1})}$ is sufficiently small.
\end{remark}

\subsection{Organization of the manuscript}

\label{Subsec-organization}The following outlines the detailed structure of
the present paper. Section \ref{section notations and preliminaries}
presents preliminary results for Morrey, Zorko, and block spaces, starting
with notations in Subsection \ref{subsection notations}. In Subsections \ref%
{subsections Embeddings and Approximation by Smooth Functions} and \ref%
{subsection Fractional and Singular Integral Operators}, we establish some
propositions on embeddings, mollification, and integral operators in these
spaces.

Next, in Section \ref{section:NewAuxiliaryResults}, we extend some classical
results in analysis and PDE theory to Zorko and block space frameworks.
First, in Subsection \ref{subsection Stein Extensions in Zorko Spaces}, we
extend the Stein extension theorem to Zorko spaces in Theorem \ref%
{Theorem:SteinExtension}. As a consequence, in Subsection \ref{subsection
Compact Embeddings and Poincare Inequalities}, we establish a version of the
Rellich-Kondrachov compact embedding theorem for Zorko spaces and
subsequently the Poincar\'{e}\ inequality (see Theorems \ref%
{Theorem:RellichKondrachovZorkoSpaces} and \ref%
{Theorem:PoincareInequalityZorkoSpaces}). We then extend both theorems to
block spaces through the duality relations in %
\eqref{equation:DualityRelations} and a Bogovskii-type result (Proposition %
\ref{proposition Bogovskii}). Both the Rellich-Kondrachov theorem and the
Poincar\'{e}\ inequality are directly used in the proof of the Helmholtz
decomposition.

In Section \ref{section main result}, we establish our main result, the
Helmholtz decomposition in Morrey, Zorko, and block spaces (Theorem \ref%
{theorem:HelmholtzDecomposition}). We begin with Zorko spaces, where
compactly supported smooth functions are dense, and then extend the proof by
duality to the other two spaces. This decomposition is obtained by solving
an equivalent weak Neumann problem (see Lemma \ref{lemma de equivalencia}).
In Subsection \ref{subsection Helmholtz Decomposition for Zorko spaces in
half space}, we employ Green (Neumann) functions to achieve the
decomposition in $\mathbb{R}^{n}$ and $\mathbb{R}_{+}^{n}$. In Subsection %
\ref{subsection Helmholtz Decomposition for Zorko spaces in bounded or
exterior domains}, we prove the decomposition for bounded or exterior $C^{1}$
domains (see Theorem \ref%
{Theorem:HelmholtzDecompositionZorkoSpacesBoundedExteriorDomain}) by
establishing an auxiliary variational inequality (see Theorem \ref%
{variational inequality theorem}). This is accomplished using localization
and flattening arguments by following \cite{simader1992new} and the
preliminary results developed in Section \ref{section notations and
preliminaries}. Finally, in Subsection \ref{subsection Helmholtz
Decomposition for Morrey and Block spaces}, we use duality arguments to
extend the Helmholtz decomposition to Morrey and block spaces.

\section{Preliminaries}

\label{section notations and preliminaries}In this section, we give some
basic notations and present some preliminary results useful for our ends.
This includes topics such as embeddings, approximation by smooth functions,
and fractional and singular integral operators in the environment we are
working in. These topics are organized into three subsections, as previously
explained in \textit{Organization of the manuscript.} As usual, we start by
introducing the basic notations.

\subsection{Basic notations}

\label{subsection notations}

Let $\Omega \subset \mathbb{R}^n$ be a non-empty open set. Given a subset $S
\subset \mathbb{R}^n$, we write $S \subset \subset \Omega$ if $S$ is
compactly contained in $\Omega$, that is, if its closure $\overline{S}$ is
compact and $\overline{S} \subset \Omega$. The space $C^\infty(\Omega)$
(resp. $C^\infty(\overline{\Omega})$) consists of all smooth functions $f$
defined on $\Omega$ (resp. $\overline{\Omega}$). We observe that $f \in
C^\infty(\overline{\Omega})$ if and only if $f$ is infinitely differentiable
on $\Omega$ and $f$ and all its derivatives have continuous extensions to $%
\overline{\Omega}$. By the Whitney extension theorem, this condition is
equivalent to $f = g|_\Omega$ for some smooth function $g$ defined on an
open neighborhood of $\overline{\Omega}$. The space $C^\infty_0(\Omega)$
(resp. $C^\infty_0(\overline{\Omega})$) consists of functions $f$ defined on
$\Omega$ (resp. $\overline{\Omega}$) with support $\mathrm{supp}(f) :=
\overline{\{x \in \Omega \mid f(x) \neq 0\}}$ compactly contained in $\Omega$
(resp. $\overline{\Omega}$). Moreover, $C^\infty_{0,\sigma}(\Omega)$ denotes
the vector fields $\mathbf{v} \in C^\infty_0(\Omega)^n$ such that $\nabla
\cdot \mathbf{v}$ (the divergence of $\mathbf{v}$) is zero.

For any two function $f, g:\Omega \to \mathbb{R}$, by $f\ast g$ we mean the
convolution of $f$ and $g$:
\begin{equation*}
f\ast g(x):=\int_{\Omega} f(x-y)g(y)dy, \quad x \in \Omega,
\end{equation*}
whenever the integral above makes sense. We denote by $\Gamma$ the
fundamental solution of the Laplacian $\Delta$ in $\mathbb{R}^n$, $n \geq 3$%
, that is,
\begin{equation*}
\Gamma(x)=\frac{|x|^{2-n}}{(2-n)\sigma_n}, \quad x \in \mathbb{R}^n,
\end{equation*}
where $\sigma_n$ is the surface area of the unit sphere $\{x \in \mathbb{R}%
^n:|x|=1\}$. Then, for functions $f$ defined in $\mathbb{R}^n$ with suitable
decay at infinity, we have $\Delta(\Gamma\ast f)=f$.

For a set $E\subset \mathbb{R}^n$, we denote by $\mathrm{diam}(E)$ its
diameter, that is,
\begin{equation*}
\mathrm{diam}(E):=\sup\{|x-y|:x,y \in E\}.
\end{equation*}
If $E$ is Lebesgue-measurable, we denote by $|E|$ its measure.

If $X$ is a complex or real normed vector space, we denote by $X^\ast$ its
dual space, formed by the continuous linear functionals defined in $X$ and
equipped with the usual norm. If $T$ is a linear operator in $X$, $T^\ast$
is its adjoint in $X^\ast$. For $S \subset X$ and $R\subset X^\ast$, the
annihilator of $S$ is the space
\begin{equation*}
S^\perp := \{\phi \in X^\ast:\phi(x)=0 \text{ for all }x \in S\}
\end{equation*}
and the preannihilator of $R$ is the space
\begin{equation*}
^\perp R:=\{x \in X:\phi(x)=0 \text{ for all }\phi \in R\}.
\end{equation*}
If $Y$ is a normed vector space continuously embedding on $X$, that is, if $%
Y\subset X$ and $\|y\|_{X}\leq c \|y\|_{Y}$ for all $y \in Y$, we write $Y
\hookrightarrow X$. If the embedding is compact, in the sense that every
bounded sequence in $Y$ has a Cauchy subsequence in $X$, then we write $%
Y\hookrightarrow\hookrightarrow X$.

\subsection{Embeddings and Approximation by Smooth Functions}

\label{subsections Embeddings and Approximation by Smooth Functions}

Let $\Omega \subset \mathbb{R}^{n}$ open and non-empty. Given a function $f$
defined over $\Omega $ and denoting by $\tilde{f}$ its zero extension to $%
\mathbb{R}^{n}$, it follows from the definition of the norm in Morrey spaces
that
\begin{equation*}
\Vert f\Vert _{M_{q,\lambda }\left( \Omega \right) }\leq \Vert \tilde{f}%
\Vert _{M_{q,\lambda }\left( \mathbb{R}^{n}\right) }.
\end{equation*}%
On the other hand, given $x\in \mathbb{R}^{n}$ and $R>0$, such that $\Omega
_{R}(x):=B(x,R)\cap \Omega \neq \empty$, let $y\in \Omega _{R}(x)$. Then
\begin{equation*}
\Vert \tilde{f}\Vert _{L^{q}(B(x,R))}\leq \Vert f\Vert _{L^{q}(\Omega
_{2R}(y))}\leq (2R)^{\lambda /q}\Vert f\Vert _{M_{q,\lambda }\left( \Omega
\right) }.
\end{equation*}%
Then
\begin{equation*}
\Vert \tilde{f}\Vert _{M_{q,\lambda }\left( \mathbb{R}^{n}\right) }\leq
2^{\lambda /q}\Vert f\Vert _{M_{q,\lambda }\left( \Omega \right) }.
\end{equation*}%
Therefore, $f\in M_{q,\lambda }\left( \Omega \right) $ if and only if $%
\tilde{f}\in M_{q,\lambda }\left( \mathbb{R}^{n}\right) $ and we can
alternatively consider $M_{q,\lambda }\left( \Omega \right) $ as the
subspace of $M_{q,\lambda }\left( \mathbb{R}^{n}\right) $ formed by the
functions that vanish outside $\Omega $. The same remark holds for block
spaces $H_{q,\lambda }(\Omega )$ and, with approximation arguments, we also
conclude it for Zorko subspaces $\mathring{M}_{q,\lambda }\left( \Omega
\right) $. This allow us to extend for general domains $\Omega $ several
results about these spaces stated in the available literature for $\mathbb{R}%
^{n}$.

\begin{definition}[Muckenhoupt Weights and Weighted $L^{q}$ Spaces]
\label{definition muckenhoupt class} Given $1<q<\infty $, a non-negative
function $w\in L_{\mathrm{loc}}^{1}(\mathbb{R}^{n})$ is called a \emph{%
Muckenhoupt $A_{q}$-weight} if and only if there is a $c\geq 0$ such that
\begin{equation}
\left( |Q|^{-1}\int_{Q}w(x)dx\right) \left( |Q|^{-1}\int_{Q}w(x)^{-\frac{1}{%
q-1}}dx\right) ^{q-1}\leq c  \label{Muckenhoupt condition}
\end{equation}%
for all cubes $Q\subset \mathbb{R}^{n}$. For $\Omega \subset \mathbb{R}^{n}$
open and non-empty, the \emph{weighted $L^{q}(\Omega )$ space} with weight $%
w\in A_{q}$ is defined by
\begin{equation*}
L_{w}^{q}(\Omega ):=\{f\in L_{\mathrm{loc}}^{1}(\overline{\Omega }%
)\,;\,\Vert f\Vert _{L_{w}^{q}}^{q}:=\int_{\Omega }|f|^{q}wdx<\infty \}
\end{equation*}
\end{definition}

In what follows, we collect some useful continuous embeddings involving
Morrey, Zorko, block and weighted $L^{q}$ spaces.

\begin{proposition}[Embeddings]
\label{embeddings from Kato} Let $\Omega \subset \mathbb{R}^{n}$ open and
non-empty, and let $1<q_{0}\leq q_{1}<\infty $ and $0\leq \lambda _{1}\leq
\lambda _{0}<\kappa <n$. such that
\begin{equation*}
\frac{n-\lambda _{0}}{q_{0}}=\frac{n-\lambda _{1}}{q_{1}}.
\end{equation*}%
Then, the following continuous embeddings hold:

\begin{itemize}
\item[(I)] $M_{q_{1},\lambda _{1}}(\Omega )\hookrightarrow M_{q_{0},\lambda
_{0}}(\Omega )\quad \text{and}\quad \mathring{M}_{q_{1},\lambda _{1}}(\Omega
)\hookrightarrow \mathring{M}_{q_{0},\lambda _{0}}(\Omega ).$

\item[(II)] $H_{q_{0}^{\prime },\lambda _{0}}(\Omega )\hookrightarrow
H_{q_{1}^{\prime },\lambda _{1}}(\Omega )$ (where the prime symbol stands
for the conjugated exponent).

\item[(III)] $L^{\frac{n}{\alpha }}(\Omega )\hookrightarrow \mathring{M}%
_{q_{0},\lambda _{0}}(\Omega ),$ with $\alpha :=(n-\lambda )/q_{0}$.

\item[(IV)] $M_{q_{0},\lambda _{0}}(\Omega )\hookrightarrow L_{w}^{q}(\Omega
),$ with $w(x):=\left( 1+|x|^{2}\right) ^{-\kappa /2}$ $\left( \in
A_{q}\right) $.

\item[(V)] $L_{w^{\prime }}^{q^{\prime }}(\Omega )\hookrightarrow
H_{q^{\prime },\lambda }\left( \Omega \right) \hookrightarrow L^{\frac{n}{%
n-\alpha }}(\Omega ),$ where $w^{\prime }:=w^{-\frac{1}{q-1}}\left( \in
A_{q^{\prime }}\right) $.
\end{itemize}
\end{proposition}

\noindent {\textbf{Proof.}} The proof for (I) and (IV) can be found in \cite%
{kato1992strong}, while (III) is the particular case of (I) when $\lambda
_{1}=0$. Then (II) and (V) follows by duality. \fin

\begin{remark}
It is known that $C_{0}^{\infty }(\mathbb{R}^{n})$ is dense in $H_{q^{\prime
},\lambda }\left( \mathbb{R}^{n}\right) $ (\cite[Theorem 345]%
{sawano2020morrey}). Then, $C_{0}^{\infty }(\overline{\Omega })$ is dense in
$H_{q^{\prime },\lambda }\left( \Omega \right) $ for any non-empty open $%
\Omega \subset \mathbb{R}^{n}$. Moreover, despite we have defined $\mathring{%
M}_{q,\lambda }\left( \Omega \right) $ as the closure of $C_{0}^{\infty }(%
\overline{\Omega })$ in $M_{q,\lambda }\left( \Omega \right) $, we note that
$C_{0}^{\infty }(\Omega )$ is dense in in $\mathring{M}_{q,\lambda }\left(
\Omega \right) $. Indeed, by Proposition \ref{embeddings from Kato}, for $%
\alpha =(n-\lambda )/q$, we have
\begin{equation*}
C_{0}^{\infty }(\overline{\Omega })\subset L^{n/\alpha }(\Omega )\subset
\overline{C_{0}^{\infty }(\Omega )}^{L^{n/\alpha }}\subset \overline{%
C_{0}^{\infty }(\Omega )}^{M_{q,\lambda }}.
\end{equation*}%
Taking the $M_{q,\lambda }\left( \Omega \right) $-closure, it follows that
\begin{equation*}
\mathring{M}_{q,\lambda }\left( \Omega \right) =\overline{C_{0}^{\infty
}(\Omega )}^{M_{q,\lambda }}.
\end{equation*}
\end{remark}

\begin{definition}[Sobolev-Morrey, Zorko and Block Spaces]
For $m\in \mathbb{N}$, we denote by $W^{m}M_{q,\lambda }\left( \Omega
\right) $ the \emph{Sobolev-Morrey space of order $m$}, that is, the space
of all functions $f\in M_{q,\lambda }\left( \Omega \right) $ whose weak
partial derivatives $\partial ^{\beta }f$ exist and belong to $M_{q,\lambda
}\left( \Omega \right) $ for all multi-indices $\beta $ with $|\beta |\leq m$%
. The space $W^{m}M_{q,\lambda }\left( \Omega \right) $ is a Banach space
equipped with the natural norm:
\begin{equation*}
\Vert f\Vert _{W^{m}M_{q,\lambda }\left( \Omega \right) }:=\max_{|\beta
|\leq m}\{\Vert \partial ^{\beta }f\Vert _{M_{q,\lambda }\left( \Omega
\right) }\}.
\end{equation*}%
We also denote by $\nabla ^{m}f$ the tensor with components $\partial
^{\beta }f$, $|\beta |=m$, and
\begin{equation*}
\Vert \nabla ^{k}f\Vert _{M_{q,\lambda }\left( \Omega \right)
}:=\max_{|\beta |=m}\{\Vert \partial ^{\beta }f\Vert _{M_{q,\lambda }\left(
\Omega \right) }\}.
\end{equation*}%
Then, we have
\begin{equation*}
\Vert f\Vert _{W^{m}M_{q,\lambda }\left( \Omega \right) }:=\max_{0\leq k\leq
m}\{\Vert \nabla ^{k}f\Vert _{M_{q,\lambda }\left( \Omega \right) }\}.
\end{equation*}%
Moreover, $W_{0}^{m}M_{q,\lambda }\left( \Omega \right) $ denotes the
closure of $C_{0}^{\infty }(\Omega )$ in $W^{m}M_{q,\lambda }\left( \Omega
\right) $. The respective Sobolev spaces based on $\mathring{M}_{q,\lambda
}\left( \Omega \right) $ or $H_{q,\lambda }(\Omega )$ are defined in the
same way. Although we have $W_{0}^{m}\mathring{M}_{q,\lambda }\left( \Omega
\right) =W_{0}^{m}M_{q,\lambda }\left( \Omega \right) $, for convenience, we
shall keep both notations.
\end{definition}

The use of mollifiers for approximation arguments is particularly effective
in Zorko spaces, as shown in the following proposition.

\begin{proposition}[Mollifications]
\label{Proposition:Mollifications} Let $\phi \in C_{0}^{\infty }(\mathbb{R}%
^{n})$ be a non-negative function such that $\phi (x)=0$ for $|x|\geq 1$ and
$\int \phi =1$, and consider the \emph{mollifier} $\phi _{\epsilon
}(x):=\epsilon ^{-n}\phi (x/\epsilon )$, $\epsilon >0$.

\begin{itemize}
\item If $f\in \mathring{M}_{q,\lambda }\left( \Omega \right) $, then $\Vert
\phi _{\epsilon }\ast f-f\Vert _{\mathring{M}_{q,\lambda }\left( \Omega
\right) }\rightarrow 0,\quad \text{as }\epsilon \rightarrow 0.$

\item If $f\in W^{m}\mathring{M}_{q,\lambda }\left( \Omega \right) $, then $%
\Vert \phi _{\epsilon }\ast f-f\Vert _{W^{m}\mathring{M}_{q,\lambda }\left(
\Omega ^{\prime }\right) }\rightarrow 0,\quad \text{as }\epsilon \rightarrow
0,$ for all $\Omega ^{\prime }\subset \subset \Omega $.
\end{itemize}

The same properties hold if we replace $\mathring{M}_{q,\lambda }$ by $%
H_{q,\lambda }$.
\end{proposition}

\noindent{\textbf{Proof.}} We prove only for $\mathring{M}%
_{q,\lambda}\left(\Omega\right)$ spaces, since the proof for block spaces is
analogous.

Let $f \in \mathring{M}_{q,\lambda}\left(\Omega\right)$. Given $\delta>0$
arbitrarily small, let $\psi \in C^\infty_0(\overline{\Omega})$ such that $%
\|\psi - f \|_{\mathring{M}_{q,\lambda}\left(\mathbb{R}^n\right)}<\delta$.
Since $\|\phi_\epsilon\|_{L^1(\Omega)}=1$ and by H\"{o}lder inequality, we
have
\begin{align*}
|\phi_\epsilon \ast (\psi - f)(x)|^q&= \left|\int_{\mathbb{R}%
^n}\phi_\epsilon(x-y)^{1/q^{\prime}}\phi_\epsilon(x-y)^{1/q}\left(%
\psi(y)-f(y)\right)dy\right|^q \\
&\leq \int_{\mathbb{R}^n}\phi_\epsilon(x-y)\left|\psi(y)-f(y)\right|^q dy \\
& =\phi_\epsilon \ast |\psi - f|^q(x).
\end{align*}
Then, for $x_0 \in \Omega$ and $R>0$,
\begin{align*}
\int_{\Omega_R(x_0)} |\phi_\epsilon \ast (\psi - f)(x)|^q dx & =
\int_{B(x_0,R)} \int_{\mathbb{R}^n} \phi_\epsilon(y)|\psi(x-y)-f(x-y)|^q dy
dx \\
&= \int_{\mathbb{R}^n}
\phi_\epsilon(y)\left(\int_{B(x_0,R)}|\psi(x-y)-f(x-y)|^q dx\right) dy \\
& = \int_{\mathbb{R}^n} \phi_\epsilon(y)\|\psi-f\|^q_{L^q(B(x_0+y,R))} dy \\
& \leq \int_{\mathbb{R}^n} \phi_\epsilon(y)R^\lambda\|\psi-f\|^q_{\mathring{M%
}_{q,\lambda}\left(\mathbb{R}^n\right)}dy \\
& \leq R^\lambda\|\psi-f\|^q_{\mathring{M}_{q,\lambda}\left(\mathbb{R}%
^n\right)}.
\end{align*}
Then
\begin{equation*}
\|\phi_\epsilon \ast \psi - \phi_\epsilon \ast f\|_{\mathring{M}%
_{q,\lambda}\left(\Omega\right)} \leq \|\psi-f\|_{\mathring{M}%
_{q,\lambda}\left(R^n\right)}<\delta.
\end{equation*}
Furthermore, with $\alpha=(n-\lambda)/q$, by Proposition \ref{embeddings
from Kato}, we have
\begin{equation*}
\|\phi_\epsilon \ast \psi - \psi\|_{\mathring{M}_{q,\lambda}\left(\Omega%
\right)} \leq c \|\phi_\epsilon \ast \psi - \psi\|_{L^{n/\alpha}(\Omega)}\to
0,
\end{equation*}
as $\epsilon \to 0$ by the properties of mollification in Lebesgue spaces.
Then,
\begin{equation*}
\limsup_{\epsilon\to 0}\|\phi_\epsilon \ast f - f\|_{\mathring{M}%
_{q,\lambda}\left(\Omega\right)} \leq 2\delta,
\end{equation*}
proving the approximation in $\mathring{M}_{q,\lambda}\left(\Omega\right)$.

Now for $m>0$ and $\Omega^{\prime}\subset\subset\Omega$, we just note that
\begin{align*}
\partial^\beta(\phi_\epsilon\ast f) = \phi_\epsilon \ast \partial^\beta f,
\end{align*}
for $|\beta|\leq m$ and $\epsilon<\mathrm{dist}(\Omega^{\prime},\partial
\Omega)$. Then, the claim follows from the previous case. \fin

\begin{remark}
\label{remark L infinite not dense} For $\lambda >0$, Proposition \ref%
{Proposition:Mollifications} does not work in $M_{q,\lambda }\left( \Omega
\right) $, since $\phi _{\epsilon }\ast u\in L^{\infty }(\Omega )$ for $u\in
M_{q,\lambda }\left( \Omega \right) $ and $L^{\infty }(\Omega )$ is not
dense in $M_{q,\lambda }\left( \Omega \right) $. Indeed, given $x_{0}\in
\mathbb{R}^{n}$ and $R>0$ such that $B(x_{0},R)\subset \Omega $, let $\chi
_{B(x_{0},R)}$ the characteristic function of $B(x_{0},R)$ and let $\alpha
=(n-\lambda )/q$. Then, it is not hard to verify that the function $%
|x-x_{0}|^{-\alpha }\chi _{B(x_{0},R)}(x)$ belongs to $M_{q,\lambda }\left(
\Omega \right) $ but it can not be approached by $L^{\infty }$ functions. It
is worth mentioning, however, that $C^\infty_0(\Omega) $ is a dense subset
of $M_{q,\lambda}\left(\Omega\right) $ with respect to the weak-star
topology induced by the duality relation \eqref{equation:DualityRelations}.
\end{remark}

By applying the previous proposition, we can obtain a result of
approximation by smooth functions in the frameworks of Sobolev-Zorko and
Sobolev-block spaces.

\begin{proposition}[Approximation by Smooth Functions]
\label{Proposition:ApproximationbySmoothFunctions} Let $\Omega \subset
\mathbb{R}^{n}$ be a non-empty subset, $m\geq 0$ an integer, $1<q<\infty $,
and $0\leq \lambda <n$. Then, $C^{\infty }(\Omega )\cap W^{m}\mathring{M}%
_{q,\lambda }\left( \Omega \right) $ is dense in $W^{m}\mathring{M}%
_{q,\lambda }\left( \Omega \right) $, and $C^{\infty }(\Omega )\cap
W^{m}H_{q,\lambda }(\Omega )$ is dense in $W^{m}H_{q,\lambda }(\Omega )$.
\end{proposition}

\noindent{\textbf{Proof.}} The proof is strictly similar to the proof in
\cite{adams2003sobolev}, Theorem 3.17, for Sobolev spaces, and we consider
only the case of Zorko subspaces, since it is analogous for block spaces.

Let $f \in W^m\mathring{M}_{q,\lambda}\left(\Omega\right)$. We define $%
\Omega_j=\{x \in \Omega\,:\, |x|< j \text{ and } \mathrm{dist}%
(x,\partial\Omega)>1/j\}$ and $U_j:=\Omega_{j+2}\backslash \overline{%
\Omega_{j}}$, for $j=1,2,...$ Then $\{U_j\}$ is an open cover for $\Omega$.
Let $\mathcal{C}$ be a $C^\infty$-partition of unity for $\Omega$
subordinated to $\{U_j\}$, that is, a collection of functions $\psi \in
C^\infty_0(\mathbb{R}^n)$ satisfying:

\begin{itemize}
\item For each $\psi \in \mathcal{C}$ and each $x \in \mathbb{R}^n$, we have
$0 \leq \psi(x) \leq 1$.

\item If $K$ is a compact subset of $\Omega$, then all but finitely many $%
\psi \in \mathcal{C}$ are identically zero on $K$.

\item For every $\psi \in \mathcal{C}$, we have $\mathrm{supp}(\psi) \subset
U_j$ for some index $j \in \mathbb{N}$.

\item $\sum_{\psi \in \mathcal{C}} \psi(x) = 1$ for all $x \in \Omega$.
\end{itemize}

Let $\mathcal{C}_k$ denote the finite collection of functions $\psi \in
\mathcal{C}$ such that $\mathrm{supp}(\psi)\subset U_k$ and let $\psi_k$ be
the sum of all $\psi \in \mathcal{C}_k \backslash \cup _{j=1}^{k-1} \mathcal{%
C}_j$. Then, $\psi_k \in C^\infty_0(U_k)$ and $\sum_1^\infty \psi_k = 1$ on $%
\Omega$. Further, for $\epsilon(k)>0$ small enough, we have $\mathrm{supp}%
(\phi_{\epsilon(k)}\ast(\psi_k f))\subset\Omega_{k+3}\backslash \overline{%
\Omega_{k-1}},$ $k=1,2,3,,...$, where $\phi_\epsilon$ is the mollifier
defined as in Proposition \ref{Proposition:Mollifications} and, by
convenience, $\Omega_0:=\emptyset$. Therefore, $\mathrm{supp}%
(\phi_{\epsilon(k)}\ast(\psi_k f)-\psi_k f) \subset \Omega_{k-3}\backslash%
\overline{\Omega_{k-1}}$ so, fixed $\delta>0$ arbitrarily small, by
Proposition \ref{Proposition:Mollifications}, we also can suppose that
\begin{equation*}
\|\phi_{\epsilon(k)}\ast (\psi_k f) - \psi_k f\|_{W^m\mathring{M}%
_{q,\lambda}\left(\Omega\right)}<2^{-k}\delta.
\end{equation*}
Let $g:=\sum_{k=1}^\infty \phi_{\epsilon(k)}\ast (\psi_k f)$. For any $%
\Omega^{\prime}$ compactly contained in $\Omega$, there are only a finite
number of nonzero terms in the sum. Thus $g \in C^\infty(\Omega)$. Moreover,
\begin{equation*}
\|f-g\|_{W^m\mathring{M}_{q,\lambda}\left(\Omega\right)} \leq
\sum_{k=1}^\infty \|\psi_k f- \phi_{\epsilon(k)}\ast (\psi_k f)\|_{W^m%
\mathring{M}_{q,\lambda}\left(\Omega\right)} < \delta.
\end{equation*}
\fin

\subsection{Fractional and Singular Integral Operators}

\label{subsection Fractional and Singular Integral Operators}

In this subsection, we present results on integral operators acting on
Morrey and block spaces. Although this type of result is relatively
well-known, we were unable to locate complete statements that fully align
with the specific goals of our study. Then, for the reader convenience, the
statements and proofs are provided in the next two propositions, which are
slight adaptations of Lemmas 4.1 and 4.2 in \cite{kato1992strong}.

\begin{proposition}[Fractional Integral Operators]
\label{proposition fractional integral operators} Let $\Omega$ a non-empty,
open, and bounded subset of $\mathbb{R}^n$, $0<\delta\leq n$ and consider
the operator
\begin{equation*}
(I_\delta f)(x):=\int_{\Omega}\frac{f(y)dy}{|x-y|^{n-\delta}}, \quad x \in
\Omega.
\end{equation*}
If $1< q_0,q_1< \infty$ and $0\leq \lambda_0,\lambda_1<n$ satisfy

\begin{equation}  \label{eq q}
\frac{n}{q_0}-\frac{n}{q_1}\leq \delta
\end{equation}
and
\begin{equation*}
\frac{n-\lambda_0}{q_0}-\frac{n-\lambda_1}{q_1}\leq \delta,
\end{equation*}
then $I_\delta$ is bounded from $M_{q_0,\lambda_0}(\Omega)$ to $%
M_{q_1,\lambda_1}(\Omega)$ and from $H_{q_1^{\prime},\lambda_1}(\Omega)$ to $%
H_{q_0^{\prime},\lambda_0}(\Omega)$.
\end{proposition}

\noindent{\textbf{Proof.}} We only have to prove the continuity in Morrey
spaces since the proof for block spaces follows by duality.

Initially, let us suppose that
\begin{equation}
0<\delta <\frac{n-\lambda _{0}}{q_{0}}.  \label{eq q1}
\end{equation}%
Let $f\in M_{q_{0},\lambda _{0}}(\Omega )$, $x_{0}\in \Omega $, $0<R\leq
\mathrm{diam}(\Omega )$, and $x\in \Omega _{R}(x_{0}):=\Omega \cap
B_{R}(x_{0})$. Then,
\begin{equation*}
(I_{\delta }f)(x)=\int_{\Omega _{R}(x)}\frac{f(y)dy}{|x-y|^{n-\delta }}%
+\int_{\Omega ^{R}(x)}\frac{f(y)dy}{|x-y|^{n-\delta }}:=(I_{\delta }^{\prime
}f)(x)+(I_{\delta }^{\prime \prime }f)(x),
\end{equation*}%
where $\Omega ^{R}(x):=\Omega \backslash B_{R}(x)$. By \eqref{eq q1}, we
have $\lambda _{0}/q_{0}+n/q_{0}^{\prime }<n-\delta .$ Then, there are $%
r>\lambda _{0}$ and $s>n$ such that
\begin{equation*}
r/q_{0}+s/q_{0}^{\prime }=n-\delta .
\end{equation*}%
Then,
\begin{gather*}
|(I_{\delta }^{\prime \prime }f)(x)|\leq \int_{\Omega ^{R}(x)}\frac{|f(y)|dy%
}{|x-y|^{n-\delta }}\leq c\left( \int_{\Omega ^{R}(x)}\frac{|f(y)|^{q_{0}}dy%
}{|x-y|^{r}}\right) ^{1/q_{0}}\left( \int_{\Omega
^{R}(x)}|x-y|^{-s}dy\right) ^{1/q_{0}^{\prime }} \\
\leq c\left( \int_{|y|>R}|y|^{-r}|f(x-y)|^{q_{0}}dy\right) ^{1/q_{0}}\left(
\int_{R}^{\infty }t^{-s+n-1}dt\right) ^{1/q_{0}^{\prime }}\leq c\left(
\int_{R}^{\infty }t^{-r}d\rho (t)\right) ^{1/q_{0}}R^{-\frac{s-n}{%
q_{0}^{\prime }}},
\end{gather*}%
where $\rho (t):=\Vert f\Vert _{L^{q_{0}}(B_{t}(x))}^{q_{0}}\leq t^{\lambda
_{0}}\Vert f\Vert _{M_{q_{0},\lambda _{0}}(\Omega )}^{q_{0}}$. By
integration by parts,
\begin{gather*}
|(I_{\delta }^{\prime \prime }f)(x)|\leq c\left( -\int_{R}^{\infty
}t^{-r-1}\rho (t)dt\right) ^{1/q_{0}}R^{-\frac{s-n}{q_{0}^{\prime }}} \\
\leq cR^{-\frac{r-\lambda _{0}}{q_{0}}}\Vert f\Vert _{M_{q_{0},\lambda
_{0}}(\Omega )}R^{-\frac{s-n}{q_{0}^{\prime }}}=cR^{\delta -\frac{n-\lambda
_{0}}{q_{0}}}\Vert f\Vert _{M_{q_{0},\lambda _{0}}(\Omega )}
\end{gather*}%
hence
\begin{equation}
R^{-\lambda _{1}/q_{1}}\Vert I_{\delta }^{\prime \prime }f\Vert
_{L^{q_{1}}(\Omega _{R}(x_{0}))}\leq cR^{\delta -\frac{n-\lambda _{0}}{q_{0}}%
+\frac{n-\lambda _{1}}{q_{1}}}\Vert f\Vert _{M_{\lambda _{0},q_{0}}(\Omega
)}\leq c\Vert f\Vert _{M_{\lambda _{0},q_{0}}(\Omega )}.
\end{equation}%
As for $I_{\delta }^{\prime }f$, since $\Omega _{R}(x)\subset \Omega
_{2R(x_{0})}$, we have
\begin{equation}  \label{eq 2}
|(I_{\delta }^{\prime }f)(x)|\leq \int_{\Omega _{2R}(x_{0})}\frac{|f(y)|dy}{%
|x-y|^{n-\delta }}
\end{equation}%
and so
\begin{equation*}
\Vert I_{\delta }^{\prime }f\Vert _{L^{q_{1}}(\Omega _{R}(x_{0}))}\leq
\left\Vert \int_{\Omega _{2R}(x_{0})}\frac{|f(y)|dy}{|\cdot -y|^{n-\delta }}%
\right\Vert _{L^{q_{1}}(\Omega _{R}(x_{0}))}\leq \left\Vert \int_{\Omega
_{2R}(x_{0})}\frac{|f(y)|dy}{|\cdot -y|^{n-\delta }}\right\Vert
_{L^{q_{1}}(\Omega _{2R}(x_{0}))}.
\end{equation*}

Therefore (see \cite{gilbarg1977elliptic}, Lemma 7.12),
\begin{equation*}
\|I^{\prime}_\delta f\|_{L^{q_1}(\Omega_R(x_0))}\leq c R^{\delta - \frac{n}{%
q_0}+\frac{n}{q_1}}\|f\|_{L^{q_0}(\Omega_{2R}(x_0))},
\end{equation*}
which implies
\begin{equation}  \label{eq 3}
R^{-\lambda_1/q_1}\|I^{\prime}_\delta f\|_{L^{q_1}(\Omega_R(x_0))}\leq c
R^{\delta - \frac{n-\lambda_0}{q_0}+\frac{n-\lambda_1}{q_1}%
}\|f\|_{M_{\lambda_0,q_0}(\Omega)}\leq c \|f\|_{M_{\lambda_0,q_0}(\Omega)}.
\end{equation}
From \eqref{eq 2} and \eqref{eq 3}, it follows that
\begin{equation*}
\|I_\delta f\|_{M_{\lambda_1,q_1}(\Omega)} \leq c
\|f\|_{M_{\lambda_0,q_0}(\Omega)}.
\end{equation*}

Now, let us consider the case in which
\begin{equation*}
\delta \geq \frac{n-\lambda_0}{q_0}.
\end{equation*}
If $\delta \geq n/q_0$, then (\cite{gilbarg1977elliptic}, Lemma 7.12) $%
I_{\delta}:L^{q_0}(\Omega) \to L^{\infty}(\Omega)$. Since by Proposition \ref%
{embeddings from Kato} it follows that $M_{q_0,\lambda_0}(\Omega)%
\hookrightarrow L^{q_0}(\Omega)$ and $L^{\infty}(\Omega)\hookrightarrow
M_{q_1,\lambda_1}(\Omega)$ for bounded domains, we have nothing to do. Let
us suppose
\begin{equation*}
\delta <\frac{n}{q_0}.
\end{equation*}
Then,
\begin{equation*}
\frac{n-\lambda_0}{q_0}-\frac{n-\lambda_1}{q_1}<\frac{n-\lambda_0}{q_0}\leq
\delta < \frac{n}{q_0}.
\end{equation*}
Therefore, there is $0\leq \mu <\lambda_0$ such that
\begin{equation*}
\frac{n-\mu}{q_0}-\frac{n-\lambda_1}{q_1}\leq \delta<\frac{n-\mu}{q_0}.
\end{equation*}
By the first case, $I_\delta : M_{q_0,\mu}(\Omega) \to
M_{q_1,\lambda_1}(\Omega)$ hence the proof follows since $%
M_{q_0,\lambda_0}(\Omega)\hookrightarrow M_{q_0,\mu}(\Omega)$ by Proposition %
\ref{embeddings from Kato}. \fin

\begin{proposition}[Singular Integral Operators]
\label{Proposition Calderon-Zygmund}

Let $S:\Omega\times \mathbb{R}^n \backslash \{0\}\to \mathbb{R}$ be a
singular kernel of Calder\'{o}n-Zygmund type, i.e., $S(x,y)=\nu(x,y)/|y|^n$
with
\begin{equation*}
\nu(x,y)=\nu(x,\alpha y), \quad \alpha>0,\, x \in \Omega, y \in \mathbb{R}%
^n\backslash{0},
\end{equation*}
\begin{equation*}
\int_{|y|=1}\nu(x,y)dy=0, \quad x \in \Omega,
\end{equation*}
\begin{equation*}
|\nu(x,y)|\leq C, \quad x \in \Omega, \, |y|=1.
\end{equation*}
Then, the operator $f\mapsto Tf$, where
\begin{equation*}
(Tf)(x)= \int_\Omega S(x,x-y)f(y)dy,
\end{equation*}
is bounded in $M_{q,\lambda}(\Omega)$.
\end{proposition}

\noindent{\textbf{Proof.}} Given $f \in M_{q,\lambda}(\Omega)$, $r>0$, $x_0
\in \mathbb{R}^n$, we aim to show that
\begin{equation*}
r^{-\lambda/q}\|Tf\|_{L^q(\Omega_r(x_0))}\leq c
\|f\|_{M_{q,\lambda}\left(\Omega\right)},
\end{equation*}
where $c$ does not depend upon $f, r, x_0$. For this purpose, let us denote $%
S_r(x,y):=S(x,y)$, if $|y|<r$, and $:=0$ otherwise, $%
S^r(x,y)=S(x,y)-S_r(x,y) $,
\begin{equation*}
g_r(x):=\int_\Omega S_r(x,x-y)f(y)dy, \quad \text{and}\quad
g^r(x):=\int_\Omega S^r(x,x-y)f(y)dy.
\end{equation*}
Then, $Tf=g_r+g^r$. Now, let $s_1>\lambda$ and $s_2>n$ such that $%
s_1/q+s_2/q^{\prime}=n$. We have
\begin{align*}
|g^r(x)| \leq & \int_{\mathbb{R}^n} |S^r(x,x-y)||f(y)| dy \\
= & \int_{|y|\geq r} |S(x,y)||f(x-y)| dy \\
\leq &\, C \int_{|y|\geq r} |y|^{-n} |f(x-y)| dy \\
= &\, C\int_{|y|\geq r}
|y|^{-s_2/q^{\prime}}\left(|y|^{-s_1/q}|f(x-y)|\right) dy \\
\leq &\, C \left(\int_{|y|\geq r}
|y|^{-s_2}dy\right)^{1/q^{\prime}}\left(\int_{|y|\geq
r}|y|^{-s_1}|f(x-y)|^qdy\right)^{1/q} \\
= & \, c r^{-(s_2-n)/q^{\prime}}\left(\int_{|y|\geq
r}|y|^{-s_1}|f(x-y)|^qdy\right)^{1/q} \\
= &\, c r^{-(s_2-n)/q^{\prime}} \left(\int_r^\infty
t^{-s_1}d\rho(t)\right)^{1/q},
\end{align*}
where $\rho(t)=\int_{B(x,t)}|f|^q \leq
t^{\lambda}\|f\|_{M_{q,\lambda}\left(\Omega\right)}^q$. Therefore,
\begin{align*}
|g^r(x)| \leq & \, c r^{-(s_2-n)/q^{\prime}}\left(
\left.t^{-s_1}\rho(t)\right|_{t=r}^{t=\infty} - \int_r^\infty
(-s_1)t^{-1-s_1}\rho(t)dt\right)^{1/q} \\
= & \, c r^{-(s_2-n)/q^{\prime}} \left( -r^{-s_1}\rho(r)+s_1\int_r^\infty
t^{-1-s_1}\rho(t)dt \right)^{1/q} \\
\leq & \, c r^{-(s_2-n)/q^{\prime}} \left(\int_r^\infty
t^{-1-(s_1-\lambda)}dt\right)^{1/q}\|f\|_{M_{q,\lambda}\left(\Omega\right)}
\\
= & \, cr^{-(s_2-n)/q^{\prime}-(s_1 - \lambda)/q}
\|f\|_{M_{q,\lambda}\left(\Omega\right)} \\
= & \, cr^{-(n-\lambda)/q}\|f\|_{M_{q,\lambda}\left(\Omega\right)},
\end{align*}
which implies
\begin{equation}  \label{eq qwe1}
\|g^r\|_{L^q(\Omega_r(x_0))} \leq
cr^{-(n-\lambda)/q}\|f\|_{M_{q,\lambda}\left(\Omega\right)}
|\Omega_r(x_0)|^{1/q} \leq c
r^{\lambda/q}\|f\|_{M_{q,\lambda}\left(\Omega\right)}.
\end{equation}
As for $g_r$, from the definition of $S_r$ it follows that, for $x \in
\Omega_r(x_0)$,
\begin{equation*}
g_r(x)=\int_{\mathbb{R}^n}S_{r}(x,x-y)\widetilde{f}(y)dy,
\end{equation*}
where $\widetilde{f}(y):=f(y)$ for $y \in B_{2r}(x)$ and $:=0$ otherwise.
Then, from the boundedness of Calder\'{o}n-Zygmund (maximal) singular
integral operators in $L^q(\mathbb{R}^n)$ (\cite{calderon1956singular},
Theorems 1 and 2), it follows that
\begin{equation}  \label{eq qwe2}
\|g_r\|_{L^q(\Omega_r(x_0))} \leq c \|\widetilde{f}\|_{L^q(\mathbb{R}^n)}
\leq c\|f\|_{L^q(\Omega_{2r}(x_0))} \leq c
r^{\lambda/q}\|f\|_{M_{q,\lambda}\left(\Omega\right)},
\end{equation}
where the constant $c$ is independent of $r$. From \eqref{eq qwe1} and %
\eqref{eq qwe2} the lemma follows. \fin

\section{Auxiliary Results}

\label{section:NewAuxiliaryResults} The aim of this section is to obtain
some key results essential for constructing the decomposition within our
framework. In Zorko spaces, we first prove a Stein extension-type theorem,
which allows us to also establish a Rellich-Kondrachov-type theorem and,
consequently, a Poincar\'{e}-type inequality. The latter two results are
extended to block spaces through duality arguments, and a Bogovskii-type
result is derived by following the approach in \cite{galdi2011introduction}.
The Rellich-Kondrachov theorem and the Poincar\'{e}\ inequality, in both the
Zorko and block space frameworks, play a direct role in the proof of the
Helmholtz decomposition presented in the next section.

\subsection{Stein Extensions in Zorko Spaces}

\label{subsection Stein Extensions in Zorko Spaces}This section is devoted
to extending a known result on Stein extensions (see \cite[Chapter VI,
Section 3]{stein1970singular}) to Zorko spaces. We will need to introduce
another class of subspaces of Sobolev-Morrey spaces, which we define below.

\begin{definition}[$\mathfrak{M}^{k}_{q,\protect\lambda}\left(\Omega\right)$
spaces]
\label{Definition:WSMspaces} Let $\Omega \subset \mathbb{R}^n$ be a
non-empty open set and $k\geq 0$ an integer. We denote by $\mathfrak{M}%
^{k}_{q,\lambda}\left(\Omega\right)$ the closure in $W^kM_{q,\lambda}\left(%
\Omega\right)$ of the set of all functions $f\in C^\infty(\overline{\Omega}%
)\cap W^kM_{q,\lambda}\left(\Omega\right)$ that are bounded and have all its
partial derivative bounded. If we also denote $\mathfrak{M}%
_{q,\lambda}\left(\Omega\right):=\mathfrak{M}^{0}_{q,\lambda}\left(\Omega%
\right)$. Observe that $\mathring{M}_{q,\lambda}\left(\Omega\right)%
\hookrightarrow \mathfrak{M}_{q,\lambda}\left(\Omega\right)$.
\end{definition}

\begin{definition}[Special Lipschitz Domain]
\label{Definition:SpecialLipschitzDomain} An open set $\Omega \subset
\mathbb{R}^{n}$ will be called a \emph{special Lipschitz domain} if there is
a Lipschitz function $\phi :\mathbb{R}^{n-1}\rightarrow \mathbb{R}$ such
that $\Omega =\{x=(x^{\prime },x_{n})\in \mathbb{R}^{n}\,:\,x_{n}>\phi
(x^{\prime })\}$. The Lipschitz constant of $\phi $, that is, the smallest $%
M $ such that
\begin{equation*}
|\phi (x^{\prime })-\phi (y^{\prime })|\leq M|x^{\prime }-y^{\prime }|,\quad
x^{\prime },y^{\prime }\in \mathbb{R}^{n-1}
\end{equation*}%
will be called the \emph{\ Lipschitz bound} of $\Omega $. For extension, any
domain that is congruent to a special Lipschitz domain up to a rigid
movement will also be called of special Lipschitz.
\end{definition}

The following lemma establishes a density property for the Sobolev-Zorko
space $W^{k}\mathring{M}_{q,\lambda}(\Omega)$ for a special Lipschitz domain
$\Omega$. As a consequence, we obtain an embedding from this space to $%
\mathfrak{M}_{q,\lambda }^{k}\left( \Omega \right) .$

\begin{lemma}
\label{Lema:DensitySpecialLipschitzDomain} Let $\Omega \subset \mathbb{R}%
^{n} $ be a special Lipschitz domain, $k\geq 0$ an integer, $1<q<\infty $,
and $0\leq \lambda <n$. Then the space of functions in $C^{\infty }(%
\overline{\Omega })\cap W^{k}\mathring{M}_{q,\lambda }\left( \Omega \right) $
that are bounded and have all their derivatives bounded is dense in $W^{k}%
\mathring{M}_{q,\lambda }\left( \Omega \right) $. In particular,
\begin{equation*}
W^{k}\mathring{M}_{q,\lambda }\left( \Omega \right) \hookrightarrow
\mathfrak{M}_{q,\lambda }^{k}\left( \Omega \right) .
\end{equation*}
\end{lemma}

\noindent{\textbf{Proof.}} For $M$ as defined in Definition \ref%
{Definition:MinimallySmoothDomain}, let
\begin{equation*}
\mathcal{C}_M := \{x = (x^{\prime}, x_n) \in \mathbb{R}^n \mid M
|x^{\prime}| < |x_n|, \ x_n < 0\}.
\end{equation*}
Then, $\mathcal{C}_M$ is an open cone with its vertex at the origin and
oriented downward. Let $\eta \in C^\infty(\mathbb{R}^n)$ be a non-negative
function such that $\int_{\mathbb{R}^n} \eta = 1$ and $\mathrm{supp}(\eta)
\subset \mathcal{C}_M$. Define $\eta_\epsilon(x) := \epsilon^{-n} \eta(x /
\epsilon)$. Since $\mathrm{supp}(\eta_\epsilon) \subset \mathcal{C}_M$, for
any $f \in \mathring{M}_{q,\lambda}\left(\Omega\right)$, the convolution $%
\eta_\epsilon \ast f$ is well-defined and smooth in an open neighborhood of $%
\overline{\Omega}$. Moreover,
\begin{equation*}
|\partial^\beta f_\epsilon(x)| = \left| \int_{\mathbb{R}^n} (\partial^\beta
\eta_\epsilon)(x-y)f(y)\,dy\right| \leq c\|\partial^\beta
\eta_\epsilon\|_{H_{q^{\prime},\lambda}\left(\mathbb{R}^n\right)}\|f\|_{%
\mathring{M}_{q,\lambda}\left(\Omega\right)} \leq c(\beta, \epsilon),
\end{equation*}
for $x \in \Omega$ and any multi-index $\beta$. Thus, $\partial^\beta
f_\epsilon$ is bounded. Finally, the same steps used in the proof of
Proposition \ref{Proposition:Mollifications} show that $\eta_\epsilon \ast f
\to f$ in $W^k \mathring{M}_{q,\lambda}\left(\Omega\right)$.

\fin

Before proceeding, we recall three lemmas, the proofs of which can be found
in \cite[Chapter VI, Section 3]{stein1970singular}.

\begin{lemma}[Regularized Distance]
\label{Lemma:RegularizedDistance} Let $F \subset \mathbb{R}^n$ be a proper
non-empty closed set and let $\delta(x)$ denote the distance from a point $x
\in \mathbb{R}^n$ to $F$. Then there exists a function $\theta(x)=%
\theta(x,F) $ defined in $\mathbb{R}^n\backslash F$ such that

\begin{itemize}
\item $c_1 \delta(x) \leq \theta (x) \leq c_2 \delta(x), \quad x \in \mathbb{%
R}^n\backslash F.$

\item $\theta$ is $C^\infty$ in $\mathbb{R}^n \backslash F$ and
\begin{equation*}
|\partial^\beta \theta(x)|\leq c_\beta \delta(x)^{1-|\beta|},
\end{equation*}
for all multi-indexes $\beta \in (\mathbb{N}_0)^n.$

The constants $c_1,c_2$, and $c_\beta$ are independent of $F$.
\end{itemize}
\end{lemma}

\begin{lemma}
\label{lema function stein} There exists a continuous function $\psi$
defined on $[1,\infty)$ such that $\psi(t)=O(t^{-N})$ as $t\to \infty$ for
every $N>0$, and which satisfies
\begin{equation*}
\int_1^\infty\psi(t)dt = 1, \quad \int_1^\infty t^k\psi(t)dt=0, \quad \text{%
for }k=1,2,...
\end{equation*}
\end{lemma}

\begin{lemma}
\label{lema constant stein} Let $\Omega \subset \mathbb{R}^{n}$ be a special
Lipschitz domain, let us say, $\Omega =\{x=(x^{\prime },x_{n})\in \mathbb{R}%
^{n}\,;\,x_{n}>\phi (x^{\prime })\}$, and let $\theta (x)=\theta (x,%
\overline{\Omega })$ be the regularized distance from $\overline{\Omega }$.
Then there exists a constant $m>0$, which depends only on the Lipschitz
bound of $\Omega $, so that $m\theta (x)\geq \phi (x^{\prime })-x_{n}$ for
all $x=(x^{\prime },x_{n})\in \mathbb{R}^{n}\backslash \overline{\Omega }.$
\end{lemma}

We are now in a position to prove a result on the Stein extension in the
context of Zorko and $\mathfrak{M}_{q,\lambda }^{k}$ spaces over a special
Lipschitz domain $\Omega $, as stated in the following lemma. Later, through
a localization argument, we extend the result to the cases of bounded and
exterior Lipschitz domains.

\begin{lemma}
\label{Lemma:SteinExtensionSpecial} Let $\Omega \subset \mathbb{R}^{n}$ be a
special Lipschitz domain, $1<q<\infty $ and $0\leq \lambda <n$. Then, there
is a bounded linear operator $E:\mathfrak{M}_{q,\lambda }\left( \Omega
\right) \rightarrow \mathfrak{M}_{q,\lambda }\left( \mathbb{R}^{n}\right) $
with the following properties.

\begin{itemize}
\item[(I)] $E$ is an extension operator, that is, $E(f)|_{\Omega }=f$ for
all $f\in \mathfrak{M}_{q,\lambda }\left( \Omega \right) $.

\item[(II)] For each $k\in \mathbb{N}_{0}$ and $f\in \mathfrak{M}_{q,\lambda
}^{k}\left( \Omega \right) $, we have
\begin{equation}
\Vert \nabla ^{k}E(f)\Vert _{M_{q,\lambda }\left( \mathbb{R}^{n}\right)
}\leq c(k)\Vert \nabla ^{k}f\Vert _{M_{q,\lambda }\left( \Omega \right) }.
\label{estimate stein lemma}
\end{equation}%
In particular, $E$ maps $\mathfrak{M}_{q,\lambda }^{k}\left( \Omega \right) $
continuously into $\mathfrak{M}_{q,\lambda }^{k}\left( \mathbb{R}^{n}\right)
$.

\item[(III)] For each $k\in \mathbb{N}_{0}$, $E$ maps $W^{k}\mathring{M}%
_{q,\lambda }\left( \Omega \right) $ continuously into $\mathfrak{M}%
_{q,\lambda }^{k}\left( \mathbb{R}^{n}\right) $.
\end{itemize}
\end{lemma}

\noindent{\textbf{Proof.}} Recall that, by Definition \ref%
{Definition:SpecialLipschitzDomain}, there exists a Lipschitz function $\phi$
with Lipschitz constant $M$ such that, up to a rigid movement, $\Omega = \{x
= (x^{\prime},x_n) \in \mathbb{R}^n \,;\, x_n > \phi(x^{\prime})\}$. Then,
for $f \in C^\infty (\overline{\Omega}) \cap
W^kM_{q,\lambda}\left(\Omega\right)$, bounded and having all partial
derivatives bounded, we define
\begin{equation}  \label{eq def ext}
E(f)(x^{\prime},x_n) := \int_1^\infty f(x^{\prime},x_n + t
\delta^\ast(x))\psi(t)\,dt, \quad x_n < \phi(x^{\prime}),
\end{equation}
and let $E(f) := f$ in $\overline{\Omega}$. Here, $\psi$ is the function
given by Lemma \ref{lema function stein}, and $\delta^\ast = 2m\theta$,
where $\theta$ is the regularized distance from $\overline{\Omega}$, and $m$
is given by Lemma \ref{lema constant stein}. Then, (I) follows immediately.
Furthermore, by Lemma \ref{Lema:DensitySpecialLipschitzDomain}, $W^k%
\mathring{M}_{q,\lambda}\left(\Omega\right)\hookrightarrow \mathfrak{M}%
^{k}_{q,\lambda}\left(\Omega\right)$. Then, (III) follows from (II). Let us
prove (II).

By Definition \ref{Definition:WSMspaces}, the space of functions $f$ as
above is dense in $\mathfrak{M}_{q,\lambda }^{k}\left( \Omega \right) $.
Then, it is enough to show \eqref{estimate stein lemma} for such a function $%
f$. We start by noting that $E(f)\in C^{\infty }(\overline{\mathbb{R}^{n}})$
and that it, together with all its partial derivatives, is bounded. The
proof of this relies on the fact that $f$ has this property in $\overline{%
\Omega }$ and does not depend on the assumption that $f\in W^{k}M_{q,\lambda
}\left( \Omega \right) $. We briefly outline the idea here and refer to the
proof of Theorem 5' in \cite{stein1970singular} for more details. The main
idea is to show that, on the boundary $\partial \Omega $, the partial
derivatives of $E(f)$ coming from $\mathbb{R}^{n}\backslash \overline{\Omega
}$ coincide with the corresponding partial derivatives of $f$ coming from $%
\Omega $. To verify this, let us consider $\partial ^{\alpha }E(f)(x)$,
where $x\in \mathbb{R}^{n}\backslash \overline{\Omega }$. The case $|\alpha
|=0$ is straightforward, so we assume $|\alpha |>0$. First, we note (this
can be checked by induction on $|\alpha |$) that the partial derivative $%
\partial ^{\alpha }\left( f(x^{\prime },x_{n}+t\delta ^{\ast }(x))\right) $
is the sum of $\left( \partial ^{\alpha }f\right) (x^{\prime },x_{n}+t\delta
^{\ast }(x))$ and a linear combination of terms of the form
\begin{equation}
(\partial ^{\beta }f)(x^{\prime },x_{n}+t\delta ^{\ast
}(x))\,t^{r}\,\partial ^{\gamma }g(x),  \label{eq derivative of Ef}
\end{equation}%
where $|\beta |+|\gamma |=|\alpha |$, $g$ is a monomial in the first-order
partial derivatives of $\delta ^{\ast }$, that is,
\begin{equation*}
g=(\partial _{1}\delta ^{\ast })^{r_{1}}(\partial _{2}\delta ^{\ast
})^{r_{2}}\cdots (\partial _{n}\delta ^{\ast })^{r_{n}},
\end{equation*}%
for some integers $r_{1},\dots ,r_{n}\geq 0$, and $r=r_{1}+\cdots +r_{n}>0$.
Then, by Lemma \ref{lema function stein}, the difference
\begin{equation*}
\partial ^{\alpha }E(f)(x)-(\partial ^{\alpha }f)(x^{\prime },x_{n}+\delta
^{\ast }(x))
\end{equation*}%
is a linear combination of terms of the form
\begin{equation}
\int_{1}^{\infty }(\partial ^{\beta }f)(x^{\prime },x_{n}+t\delta ^{\ast
}(x))\,t^{r}\,\psi (t)\,dt\,\partial ^{\gamma }g(x).
\label{eq derivative int}
\end{equation}%
We observe that, by Lemma \ref{Lemma:RegularizedDistance}, we have
\begin{equation}
|\partial ^{\gamma }g|\leq c\,\delta ^{-|\gamma |}.  \label{estimate g}
\end{equation}%
Moreover, by Taylor's theorem with an integral remainder,
\begin{align*}
(\partial ^{\beta }f)(x^{\prime },x_{n}+t\delta ^{\ast }(x))=&
\sum_{j=0}^{|\gamma |}\frac{(t\delta ^{\ast }(x)-\delta ^{\ast }(x))^{j}}{j!}%
(\partial _{n}^{j}\partial ^{\beta }f)(x^{\prime },x_{n}+\delta ^{\ast }(x))
\\
& +\int_{\delta ^{\ast }(x)}^{t\delta ^{\ast }(x)}\frac{(s-\delta ^{\ast
}(x))^{|\gamma |}}{|\gamma |!}(\partial _{n}^{|\gamma |+1}\partial ^{\beta
}f)(x^{\prime },x_{n}+s)\,ds.
\end{align*}%
Substituting this expression into \eqref{eq derivative int}, by the
orthogonality relations from Lemma \ref{lema function stein}, all terms
vanish except for the corresponding to the remainder:
\begin{equation*}
\int_{1}^{\infty }\int_{\delta ^{\ast }(x)}^{t\delta ^{\ast }(x)}\frac{%
(s-\delta ^{\ast }(x))^{|\gamma |}}{|\gamma |!}(\partial _{n}^{|\gamma
|+1}\partial ^{\beta }f)(x^{\prime },x_{n}+s)\,ds\,t^{r}\,\psi
(t)\,dt\,\partial ^{\gamma }g(x),
\end{equation*}%
which can be estimated by
\begin{equation*}
c\Vert \nabla ^{|\alpha |+1}f\Vert _{L^{\infty }(\Omega )}\int_{1}^{\infty
}(t-1)^{|\gamma |+1}t^{r}|\psi (t)|\,dt\,\delta ^{\ast }(x)^{|\gamma
|+1}|\partial ^{\gamma }g(x)|.
\end{equation*}%
The above integral converges since $\psi $ decays sufficiently fast.
Moreover, by \eqref{estimate g} and Lemma \ref{Lemma:RegularizedDistance},
we have $\delta ^{\ast }(x)^{|\gamma |+1}|\partial ^{\gamma }g(x)|\leq
c\delta (x)\rightarrow 0$. Therefore,
\begin{equation*}
\partial ^{\alpha }E(f)(x)-(\partial ^{\alpha }f)(x^{\prime },x_{n}+\delta
^{\ast }(x))\rightarrow 0,
\end{equation*}%
as $x\rightarrow x_{0}$, meaning that $\partial ^{\alpha }E(f)(x)\rightarrow
(\partial ^{\alpha }f)(x_{0})$.

Next, we aim to show that, for $k\in \mathbb{N}_{0}$,
\begin{equation}
\Vert \nabla ^{k}E(f)\Vert _{M_{q,\lambda }\left( \mathbb{R}^{n}\right)
}\leq c(k)\Vert \nabla ^{k}f\Vert _{M_{q,\lambda }\left( \Omega \right) }.
\label{equation:InequalityinSteinLemma}
\end{equation}%
Note that this implies $\Vert E(f)\Vert _{\mathfrak{M}_{q,\lambda
}^{k}\left( \Omega \right) }\leq c\Vert f\Vert _{\mathfrak{M}_{q,\lambda
}^{k}\left( \Omega \right) }$, thus establishing (II) by density. We begin
by considering \eqref{equation:InequalityinSteinLemma} with $k=0$. First,
observe that
\begin{equation*}
2(\phi (x^{\prime })-x_{n})\leq \delta ^{\ast }(x)\leq c\delta (x)\leq
c(\phi (x^{\prime })-x_{n}).
\end{equation*}%
Then, fixing $x\in \mathbb{R}^{n}\backslash \overline{\Omega }$ and
substituting $t=(1+s)(\phi (x^{\prime })-x_{n})/\delta ^{\ast }(x)$ into
\eqref{eq
def ext}, we obtain
\begin{align*}
|E(f)(x^{\prime },x_{n})|& \leq \int_{\frac{\delta ^{\ast }(x)}{\phi
(x^{\prime })-x_{n}}-1}^{\infty }\left\vert f(x^{\prime },(\phi (x^{\prime
})-x_{n})s+\phi (x^{\prime }))\right\vert \left\vert \psi \left( \frac{%
(1+s)(\phi (x^{\prime })-x_{n})}{\delta ^{\ast }(x)}\right) \right\vert
\frac{\phi (x^{\prime })-x_{n}}{\delta ^{\ast }(x)}\,ds \\
& \leq \,c\int_{1}^{\infty }\left\vert f(x^{\prime },(\phi (x^{\prime
})-x_{n})s+\phi (x^{\prime }))\right\vert \left( \frac{(1+s)(\phi (x^{\prime
})-x_{n})}{\delta ^{\ast }(x)}\right) ^{-N}\frac{\phi (x^{\prime })-x_{n}}{%
\delta ^{\ast }(x)}\,ds, \\
& \leq c\int_{1}^{\infty }\left\vert f(x^{\prime },(\phi (x^{\prime
})-x_{n})s+\phi (x^{\prime }))\right\vert (1+s)^{-N}ds \\
& \leq c\int_{1}^{\infty }\left\vert f(x^{\prime },(\phi (x^{\prime
})-x_{n})s+\phi (x^{\prime }))\right\vert s^{-N}ds,
\end{align*}%
where $N$ is chosen sufficiently large and $c=c(N)$ (see Lemma \ref{lema
function stein}). Now, let $R>0$ and $y\in \mathbb{R}^{n}$, and consider the
cube $\tilde{Q}=[y_{1}-R,y_{1}+R]\times ...\times \lbrack y_{n}-R,y_{n}+R]$.
For simplicity, denote $\tilde{Q}=Q\times \lbrack a,b]$, where $%
Q=[y_{1}-R,y_{1}+R]\times ...\times \lbrack y_{n-1}-R,y_{n-1}+R]$ and $%
[a,b]=[y_{n}-R,y_{n}+R]$. Denoting by $\xi $ the characteristic function of $%
\mathbb{R}^{n}\backslash \overline{\Omega }$, we have
\begin{align*}
\left( \int_{\widetilde{Q}}\xi |E(f)|^{q}dx\right) ^{1/q}& =\left(
\int_{Q}\int_{a}^{b}\xi (x^{\prime },x_{n})|E(f)(x^{\prime
},x_{n})|^{q}dx_{n}dx^{\prime }\right) ^{1/q} \\
& \leq \,c\left( \int_{Q}\int_{a}^{b}\left( \int_{1}^{\infty }\xi (x^{\prime
},x_{n})\left\vert f(x^{\prime },(\phi (x^{\prime })-x_{n})s+\phi (x^{\prime
}))\right\vert s^{-N}ds\right) ^{q}dx_{n}dx^{\prime }\right) ^{1/q} \\
& \leq \,c\int_{1}^{\infty }\left( \int_{Q}\int_{a}^{b}\xi (x^{\prime
},x_{n})\left\vert f(x^{\prime },(\phi (x^{\prime })-x_{n})s+\phi (x^{\prime
}))\right\vert ^{q}s^{-Nq}dx_{n}dx^{\prime }\right) ^{1/q}ds.
\end{align*}%
We shall perform the change of variables $t=(\phi (x^{\prime })-x_{n})s+\phi
(x^{\prime })$, that is, $x_{n}=\phi (x^{\prime })-(t-\phi (x^{\prime }))/s$
in the integral with respect to $x_{n}$. Consider a point $(x^{\prime
},x_{n})\in \widetilde{Q}\backslash \overline{\Omega }$. Then, $x_{n}<\phi
(x^{\prime })$, i.e., $t>\phi (x^{\prime })$ so $(x^{\prime },t)\in \Omega $%
. Moreover, since $x_{n}\geq a$, we have
\begin{gather*}
t\leq \phi (x^{\prime })+t(\phi (x^{\prime })-a)\leq \phi (y^{\prime
})+M|x^{\prime }-y^{\prime }|+s(\phi (y^{\prime })+M|x^{\prime }-y^{\prime
}|-a) \\
\leq \phi (y^{\prime })+MR+s(\phi (y^{\prime })+MR-a):=A(s).
\end{gather*}%
Similarly,
\begin{equation*}
t\geq \phi (y^{\prime })-MR+s(\phi (y^{\prime })-MR-b):=B(s).
\end{equation*}%
Thus, we obtain
\begin{equation*}
\left( \int_{\widetilde{Q}}\xi |E(f)|^{q}dx\right) ^{1/q}\leq
c\int_{1}^{\infty }\left( \int_{Q}\int_{B(s)}^{A(s)}\eta (x^{\prime
},t)\left\vert f(x^{\prime },t)\right\vert ^{q}\frac{dt}{s}dx^{\prime
}\right) ^{1/q}s^{-N}ds,
\end{equation*}%
where $\eta $ is the characteristic function of $\Omega $. Note that, for
each $s\geq 1$, we have
\begin{equation*}
A(s)-B(s)=2M(1+s)R+s(b-a)\leq 4MsR+2sR\leq csR.
\end{equation*}%
Since the sides of $Q$ have length $2R$, there exists a cube $Q^{\prime
}(s)\subset \mathbb{R}^{n}$ with sides of length $csR$ containing $Q\times
\lbrack B(s),A(s)]$. Hence,
\begin{gather*}
(2R)^{-\lambda /q}\left( \int_{\widetilde{Q}}\xi |E(f)|^{q}dx\right)
^{1/q}\leq cR^{-\lambda /q}\int_{1}^{\infty }\left( \int_{Q^{\prime
}(s)}\eta (x)\left\vert f(x)\right\vert ^{q}dx\right) ^{1/q}s^{-N-1/q}ds \\
\leq c\int_{1}^{\infty }s^{\lambda /q}\Vert f\Vert _{M_{q,\lambda }\left(
\Omega \right) }s^{-N-1/q}ds=c\Vert f\Vert _{M_{q,\lambda }\left( \Omega
\right) }.
\end{gather*}%
Taking the supremum over all $y\in \mathbb{R}^{n}$ and $R>0$, we get
\begin{equation}
\Vert E(f)\Vert _{M_{q,\lambda }\left( \mathbb{R}^{n}\backslash \Omega
\right) }\leq c\Vert f\Vert _{M_{q,\lambda }\left( \Omega \right) }.
\label{stein estimate for morrey}
\end{equation}%
Thus,
\begin{equation*}
\Vert E(f)\Vert _{M_{q,\lambda }\left( \mathbb{R}^{n}\right) }\leq c\Vert
f\Vert _{M_{q,\lambda }\left( \Omega \right) }.
\end{equation*}%
It is important to note that the dependence of the constant $c$ on the
domain $\Omega $ is solely related to its Lipschitz constant $M$.

The proof of \eqref{equation:InequalityinSteinLemma} for $k>0$ follows
similarly. Indeed, given $\alpha \in (\mathbb{N}_{0})^{n}$ a multi-index
with $|\alpha |=k$, as discussed earlier, for $x\in \mathbb{R}^{n}\backslash
\overline{\Omega }$, the partial derivative $\partial ^{\alpha }E(f)(x)$ is
the sum of $\left( \partial ^{\alpha }f\right) (x^{\prime },x_{n}+t\delta
^{\ast }(x))$ and a linear combination of terms of the form
\begin{equation}
\int_{1}^{\infty }(\partial ^{\beta }f)(x^{\prime },x_{n}+t\delta ^{\ast
}(x))\,t^{r}\,\psi (t)\,dt\,\partial ^{\gamma }g(x),  \label{eq rep df}
\end{equation}%
where $|\beta |+|\gamma |=k$, $g=(\partial _{1}\delta ^{\ast
})^{r_{1}}(\partial _{2}\delta ^{\ast })^{r_{2}}\cdots (\partial _{n}\delta
^{\ast })^{r_{n}},$ and $r=r_{1}+\cdots +r_{n}>0$. The same steps used to
prove \eqref{stein estimate for morrey} show that $\left( \partial ^{\alpha
}f\right) (x^{\prime },x_{n}+t\delta ^{\ast }(x))$ has norm $\Vert \cdot
\Vert _{M_{q,\lambda }\left( \mathbb{R}^{n}\backslash \Omega \right) }$
bounded by $c\Vert \nabla ^{k}f\Vert _{M_{q,\lambda }\left( \Omega \right) }$%
. The same holds for the terms of the form \eqref{eq rep df} with $|\gamma
|=0$, since $g$ is bounded, $|\beta |=k$, and $t^{r}\psi (t)$ decays
sufficiently fast. Finally, for $|\gamma |>0$, we apply Taylor's theorem
(expanding to the $(|\gamma |-1)$-th term) to conclude that \eqref{eq rep df}
equals
\begin{align*}
& \int_{1}^{\infty }\int_{\delta ^{\ast }(x)}^{t\delta ^{\ast }(x)}\frac{%
(s-\delta ^{\ast }(x))^{|\gamma |-1}}{(|\gamma |-1)!}(\partial _{n}^{|\gamma
|}\partial ^{\beta }f)(x^{\prime },x_{n}+s)\,ds\,t^{r}\,\psi
(t)\,dt\,\partial ^{\gamma }g(x) \\
& =c\int_{1}^{\infty }\int_{1}^{t}(s-1)^{|\gamma |-1}(\partial _{n}^{|\gamma
|}\partial ^{\beta }f)(x^{\prime },x_{n}+s\delta ^{\ast
}(x))\,ds\,t^{r}\,\psi (t)\,dt\,\delta ^{\ast }(x)^{|\gamma |}\partial
^{\gamma }g(x) \\
& =c\int_{1}^{\infty }(s-1)^{|\gamma |-1}(\partial _{n}^{|\gamma |}\partial
^{\beta }f)(x^{\prime },x_{n}+s\delta ^{\ast }(x))\left( \int_{s}^{\infty
}\,t^{r}\,\psi (t)\,dt\right) ds\,\delta ^{\ast }(x)^{|\gamma |}\partial
^{\gamma }g(x),
\end{align*}%
which can be estimated by
\begin{equation*}
c\int_{1}^{\infty }(\partial _{n}^{|\gamma |}\partial ^{\beta }f)(x^{\prime
},x_{n}+s\delta ^{\ast }(x))s^{-N}ds,
\end{equation*}%
for $N$ arbitrarily large. Thus, the same steps used to prove
\eqref{stein
estimate for morrey} apply again, leading to %
\eqref{equation:InequalityinSteinLemma}. Therefore, (II) follows by density. %
\fin

\begin{definition}[Minimally Smooth Domain]
\label{Definition:MinimallySmoothDomain} A non-empty open set $\Omega
\subset \mathbb{R}^n$ is called \emph{minimally smooth} if there exist $%
\epsilon > 0$, an integer $m > 0$, a constant $M > 0$, and a countable
(possibly finite) covering $\{U_i\}_{i \in I}$ of $\partial \Omega$,
consisting of non-empty, distinct open sets $U_i$, such that:

\begin{itemize}
\item For each $x \in \partial \Omega$, there exists some $i$ such that $%
B(x, \epsilon) \subset U_i$.

\item The intersection of any collection of more than $m$ distinct sets $U_i$
is empty.

\item For each $i$, there exists a special Lipschitz domain $\Omega_i$ with
boundary norm bounded by $M$ such that $\Omega \cap U_i = \Omega_i \cap U_i$.
\end{itemize}
\end{definition}

We note that an example of a minimally smooth domain is a bounded or
exterior Lipschitz domain, that is, where $\partial \Omega$ is locally the
graph of a Lipschitz function. In this case, $\partial \Omega$ is compact,
and the covering $\{U_i\}_{i \in I}$ in the definition above is finite.
Next, we extend Lemma \ref{Lemma:SteinExtensionSpecial} to this kind of
domain.

\begin{theorem}[Stein Extension for $\mathfrak{M}_{q,\protect\lambda%
}\left(\Omega\right)$ and $\mathring{M}_{q,\protect\lambda%
}\left(\Omega\right)$]
\label{Theorem:SteinExtension} Let $\Omega\subset \mathbb{R}^n$ be a bounded
or exterior Lipschitz domain. Then, there is a bounded linear operator $E:
\mathfrak{M}_{q,\lambda}\left(\Omega\right)\to \mathfrak{M}_{q,\lambda}\left(%
\mathbb{R}^n\right)$ with the following properties.

\begin{itemize}
\item[(I)] $E$ is a extension operator, that is, $E(f)|_{\Omega}=f$ for all $%
f \in\mathfrak{M}_{q,\lambda}\left(\Omega\right)$.

\item[(II)] For each $k \in \mathbb{N}_0$, $E$ maps both $\mathfrak{M}%
^{k}_{q,\lambda}\left(\Omega\right)$ and $W^k\mathring{M}_{q,\lambda}\left(%
\Omega\right)$ continuously into $\mathfrak{M}^{k}_{q,\lambda}\left(\mathbb{R%
}^n\right)$.
\end{itemize}
\end{theorem}

\noindent{\textbf{Proof.}} According to Definition \ref%
{Definition:MinimallySmoothDomain}, there exists a finite open covering $%
\{U_i\}_{1 \leq i \leq N}$ of $\partial \Omega$, such that for each $i$, we
have $\Omega \cap U_i = \Omega_i \cap U_i$ for some special Lipschitz domain
$\Omega_i$. Next, we take non-negative smooth functions $\Lambda_+,
\Lambda_-, \lambda_1, \dots, \lambda_N$, which are bounded and have all
their partial derivatives bounded in $\mathbb{R}^n$, satisfying the
following properties:

\begin{itemize}
\item $\Lambda_+ + \Lambda_- = 1$ in $\Omega$;

\item $\mathrm{supp}(\Lambda_-) \subset \Omega$;

\item $\sum_{i=1}^N \lambda_i(x)^2 \geq 1$ for $x$ in an open neighborhood
of $\mathrm{supp}(\Lambda_+)$;

\item $\mathrm{supp}(\lambda_i) \subset U_i$ for $i = 1, \dots, N$.
\end{itemize}

The construction of these functions is explained in detail in \cite%
{stein1970singular}, p. 190, and is omitted here for brevity.

For each $i$, let $E^i$ be the Stein extension operator for $\Omega_i$,
guaranteed by Lemma \ref{Lemma:SteinExtensionSpecial}. Given $f \in
\mathfrak{M}_{q,\lambda}\left(\Omega\right)$, we define
\begin{equation}  \label{eq sum E}
E(f)(x) := \Lambda_+(x)\left(\frac{\sum_{i=1}^N \lambda_i(x) E^i(\lambda_i
f)(x)}{\sum_{i=1}^N \lambda_i^2(x)}\right) + \Lambda_-(x)f(x), \quad x \in
\mathbb{R}^n.
\end{equation}
Note that $E(f)|_\Omega = f$. Furthermore,
\begin{align*}
\|E(f)\|_{\mathfrak{M}_{q,\lambda}\left(\mathbb{R}^n\right)} &\leq
c\left(\sum_{i=1}^N \|E^i(\lambda_i f)\|_{\mathfrak{M}_{q,\lambda}\left(%
\mathbb{R}^n\right)} + \|f\|_{\mathfrak{M}_{q,\lambda}\left(\Omega\right)}
\right) \\
&\leq c \left(\sum_{i=1}^N \|\lambda_i f\|_{\mathfrak{M}_{q,\lambda}\left(%
\Omega_i\right)} + \|f\|_{\mathfrak{M}_{q,\lambda}\left(\Omega\right)}%
\right) \leq c \|f\|_{\mathfrak{M}_{q,\lambda}\left(\Omega\right)}.
\end{align*}
Similarly, since each $E^i$ also maps $\mathring{M}_{q,\lambda}\left(%
\Omega_i\right)$ continuously into $\mathfrak{M}_{q,\lambda}\left(\mathbb{R}%
^n\right)$, we also have
\begin{equation*}
\|E(f)\|_{\mathfrak{M}_{q,\lambda}\left(\mathbb{R}^n\right)} \leq c \|f\|_{%
\mathring{M}_{q,\lambda}\left(\Omega\right)}, \quad f \in \mathring{M}%
_{q,\lambda}\left(\Omega\right).
\end{equation*}
Finally, the estimates for $E(f)$ in $\mathfrak{M}^{k}_{q,\lambda}\left(%
\mathbb{R}^n\right)$, with $k > 0$, are obtained in the same way, taking
into account the boundedness of the functions $\Lambda_+, \Lambda_-,
\lambda_1, \dots, \lambda_N$, and their partial derivatives. \fin

\subsection{Compact Embeddings and Poincar\'{e} Inequalities}

\label{subsection Compact Embeddings and Poincare Inequalities}

This part of the text focuses on obtaining some compact embeddings and
Poincar\'{e}-type inequalities in our framework. We start by recalling the
following proposition whose proof can be found in \cite[Theorem 6.1]%
{adams2015morrey}.

\begin{proposition}
\label{prop maximal operator} Let $M$ be the Hardy-Littlewood maximal
function, that is,
\begin{equation*}
M(f)(x):=\sup_{r>0}\frac{1}{r^{n}}\int_{B(x,r)}|f|,\quad x\in \mathbb{R}^{n}.
\end{equation*}%
Then, $M$ is bounded from $M_{q,\lambda }\left( \mathbb{R}^{n}\right) $ to $%
M_{q,\lambda }\left( \mathbb{R}^{n}\right) $.
\end{proposition}

The next lemma establishes a compact embedding within the framework of Zorko
spaces, first for $\Omega $ as a ball. Building on our result on Stein
extensions, we then extend this to general bounded Lipschitz domains.

\begin{lemma}
\label{lemma compact emb} Let $R>0$ and $B_{R}:=B(0,R)=\{x\in \mathbb{R}%
^{n}:|x|<R\}$. Then the following compact embedding holds:
\begin{equation}
W_{0}^{1}\mathring{M}_{q,\lambda }\left( B_{R}\right) \hookrightarrow
\hookrightarrow \mathring{M}_{q,\lambda }\left( B_{R}\right) .
\label{compact emb1}
\end{equation}
\end{lemma}

\noindent{\textbf{Proof.}} Let $\{f_n\}_{n \in \mathbb{N}} \subset
C^\infty_0(B_R)$ with $\|f_n\|_{W^1M_{q,\lambda}\left(B_R\right)}$ bounded.
Since
\begin{equation*}
W^1_0\mathring{M}_{q,\lambda}\left(B\right) \hookrightarrow
W_0^{1,q}(B_R)\hookrightarrow\hookrightarrow L^q(B_R),
\end{equation*}
there is $f \in W^{1,q}_0(B_R)$ such that, up to a subsequence, $f_n
\rightharpoonup f$ in the weak topology of $W^{1,q}_0(B_R)$ and $f_n\to f$
in $L^q(B_R)$. Let us consider the zero extension of $f$ to $\mathbb{R}^n$
which, for the sake of simplicity, we still denote by $f$. Then, $f \in
W^{1,q}_0(\mathbb{R}^n)$. We can suppose that $f \in M_{q,\lambda}\left(%
\mathbb{R}^n\right)$. Indeed, since $f_n$ is bounded in $M_{q,\lambda}\left(%
\mathbb{R}^n\right)=H_{q^{\prime},\lambda}(\mathbb{R}^n)^\ast$, by
Banach-Alaoglu theorem there is $g \in M_{q,\lambda}(\mathbb{R}^n)$ such
that, up to subsequence, $f_n\overset{\ast}{\rightharpoonup} g$ in the
weak-star topology of $H_{q^{\prime},\lambda}(\mathbb{R}^n)^\ast$. Then, for
any $\psi \in C^\infty_0(\mathbb{R}^n)$,
\begin{equation*}
\int g\psi = \lim \int f_n \psi =\int f \psi.
\end{equation*}
Therefore $f=g \in M_{q,\lambda}\left(\mathbb{R}^n\right)$. By the same
argument, we can suppose that $\nabla f \in M_{q,\lambda}\left(\mathbb{R}%
^n\right)$.

Let $\phi \in C_{0}^{\infty }(B_{1})$ with $\phi \geq 0$, $\int \phi =1$
and, for $\epsilon >0$ arbitrarily small, let $\phi _{\epsilon
}(x):=\epsilon ^{-n}\phi (x/\epsilon )$. Then, for $x\in \mathbb{R}^{n}$,
\begin{gather*}
|\phi _{\epsilon }\ast f(x)-f(x)|=\left\vert \int_{|y|<\epsilon }\phi
_{\epsilon }(y)(f(x-y)-f(x))dy\right\vert =\left\vert \int_{|y|<\epsilon
}\phi _{\epsilon }(y)\int_{0}^{1}\nabla f(x-ty)\cdot ydtdy\right\vert \\
=\epsilon \Vert \phi \Vert _{\infty }\int_{0}^{1}\int_{|y|<\epsilon
}\epsilon ^{-n}|\nabla f(x-ty)|dydt=\epsilon
c\int_{0}^{1}\int_{|x-z|<t\epsilon }\epsilon ^{-n}|\nabla f(z)|t^{-n}dzdt \\
\leq \epsilon c\sup_{r>0}\frac{1}{r^{n}}\int_{B(x,r)}|\nabla f|=\epsilon
cM(|\nabla f|)(x),
\end{gather*}%
where $M$ denotes the Hardy-Littlewood maximal function (Proposition \ref%
{prop maximal operator}). Hence
\begin{equation*}
\Vert \phi _{\epsilon }\ast f-f\Vert _{M_{q,\lambda }\left( \mathbb{R}%
^{n}\right) }\leq \epsilon c\Vert M(|\nabla f|)\Vert _{M_{q,\lambda }\left(
\mathbb{R}^{n}\right) }\leq \epsilon c\Vert \nabla f\Vert _{M_{q,\lambda
}\left( \mathbb{R}^{n}\right) }.
\end{equation*}%
Analogously,
\begin{equation*}
\Vert \phi _{\epsilon }\ast f_{n}-f_{n}\Vert _{M_{q,\lambda }\left( \mathbb{R%
}^{n}\right) }\leq \epsilon c\Vert \nabla f_{n}\Vert _{M_{q,\lambda }\left(
\mathbb{R}^{n}\right) }.
\end{equation*}%
Then,
\begin{align*}
\Vert f_{n}-f\Vert _{M_{q,\lambda }\left( B\right) }& =\Vert f_{n}-f\Vert
_{M_{q,\lambda }\left( \mathbb{R}^{n}\right) } \\
& \leq \Vert f_{n}-\phi _{\epsilon }\ast f_{n}\Vert _{M_{q,\lambda }\left(
\mathbb{R}^{n}\right) }+\Vert \phi _{\epsilon }\ast (f_{n}-f)\Vert
_{M_{q,\lambda }\left( \mathbb{R}^{n}\right) }+\Vert \phi _{\epsilon }\ast
f-f\Vert _{M_{q,\lambda }\left( \mathbb{R}^{n}\right) } \\
& \leq c\left( \epsilon \Vert \nabla f_{n}\Vert _{M_{q,\lambda }\left(
\mathbb{R}^{n}\right) }+\Vert \phi _{\epsilon }\ast (f_{n}-f)\Vert
_{L^{\infty }(\mathbb{R}^{n})}+\epsilon \Vert \nabla f\Vert _{M_{q,\lambda
}\left( \mathbb{R}^{n}\right) }\right) \\
& \leq c\left( \epsilon \Vert \nabla f_{n}\Vert _{M_{q,\lambda }\left(
\mathbb{R}^{n}\right) }+\Vert \phi _{\epsilon }\Vert _{L^{q^{\prime }}(%
\mathbb{R}^{n})}\Vert f_{n}-f\Vert _{L^{q}(\mathbb{R}^{n})}+\epsilon \Vert
\nabla f\Vert _{M_{q,\lambda }\left( \mathbb{R}^{n}\right) }\right) .
\end{align*}%
Since $\Vert \nabla f_{n}\Vert _{M_{q,\lambda }\left( \mathbb{R}^{n}\right)
} $ is bounded, by choosing $\epsilon >0$ small and then letting $%
n\rightarrow \infty $, we conclude that $\Vert f_{n}-f\Vert _{M_{q,\lambda
}\left( B_{R}\right) }\rightarrow 0$, which proves, in particular, that $%
f\in \mathring{M}_{q,\lambda }\left( B_{R}\right) :=\overline{C_{0}^{\infty
}(B_{R})}^{M_{q,\lambda }}$. Finally, \eqref{compact emb1} follows by
density. \fin

Having established the previous lemma and our result on Stein extensions
(see Theorem \ref{Theorem:SteinExtension}), we are now positioned to prove
the following Rellich-Kondrachov-type theorem.

\begin{theorem}[Rellich-Kondrachov Theorem for Zorko Spaces]
\label{Theorem:RellichKondrachovZorkoSpaces} Let $\Omega \subset \mathbb{R}%
^n $ be a bounded Lipschitz domain. Then the following compact embedding
holds:
\begin{equation}
W^1 \mathring{M}_{q,\lambda}\left(\Omega\right)
\hookrightarrow\hookrightarrow \mathring{M}_{q,\lambda}\left(\Omega\right).
\end{equation}
\end{theorem}

\noindent{\textbf{Proof.}} Let $\{f_n\}$ be a bounded sequence in $W^1%
\mathring{M}_{q,\lambda}\left(\Omega\right)$ and, for each $n$, let $Ef_n
\in \mathfrak{M}^{1}_{q,\lambda}\left(\mathbb{R}^n\right)$ be the Stein
extension of $f_n$ assured by Theorem \ref{Theorem:SteinExtension}. Then, $%
\{Ef_n\}$ is a bounded sequence in $\mathfrak{M}^{1}_{q,\lambda}\left(%
\mathbb{R}^n\right)$. Now, let $R>0$ be large enough so $B_{R/2}\supset
\Omega$, and take $\eta \in C^\infty_0(\mathbb{R}^n)$ such that $\eta=1$ in $%
B_{R/2}$ and $\eta=0$ in $\mathbb{R}^n \backslash B_R$. Since smooth
functions are dense in $\mathfrak{M}^{1}_{q,\lambda}\left(\mathbb{R}%
^n\right) $, each $g_n:=\eta Ef_n$ can be approximated in $%
W^1M_{q,\lambda}\left(B_R\right)$ by smooth, compactly supported in $B_R$,
functions. Then, the sequence $\{g_n\}$ is bounded in $W^1_0\mathring{M}%
_{q,\lambda}\left(B_R\right)$ and, by Lemma \ref{lemma compact emb}, there
exists $g \in \mathring{M}_{q,\lambda}\left(B_R\right)$ such that $%
\|g_n-g\|_{\mathring{M}_{q,\lambda}\left(B_R\right)}\to 0$ up to a
subsequence. Let $f:=g|_{\Omega}$. Then $f \in \mathring{M}%
_{q,\lambda}\left(\Omega\right)$ and
\begin{equation*}
\|f_n-f\|_{\mathring{M}_{q,\lambda}\left(\Omega\right)} \leq \|g_n-g\|_{%
\mathring{M}_{q,\lambda}\left(B_R\right)} \to 0, \quad \text{as }n \to
\infty,
\end{equation*}
which proves the compact embedding.

\fin

As a consequence of Theorem \ref{Theorem:RellichKondrachovZorkoSpaces}, we
have the following corollary, which can be seen as a first step toward
obtaining a Poincar\'{e} inequality in the context of Zorko spaces.

\begin{corollary}
\label{corollary before poincare} Let $\Omega$ a bounded and connected
Lipschitz domain in $\mathbb{R}^n$ and $\psi:\mathring{M}_{q,\lambda}\left(%
\Omega\right)\to \mathbb{R}_+$ a continuous and absolutely homogeneous
function ($\psi(sf)=|s|\psi(f)$) such that $\psi(1)\neq 0$. Then, there is a
constant $c$ such that
\begin{equation*}
\|f\|_{\mathring{M}_{q,\lambda}\left(\Omega\right)}\leq c \left(\|\nabla f
\|_{\mathring{M}_{q,\lambda}\left(\Omega\right)}+\psi(f)\right),
\end{equation*}
for all $f \in W^1\mathring{M}_{q,\lambda}\left(\Omega\right)$
\end{corollary}

\noindent {\textbf{Proof.}} Suppose by contradiction that there is a
sequence $\{f_{k}\}\subset W^{1}\mathring{M}_{q,\lambda }\left( \Omega
\right) $ such that
\begin{equation}
\Vert f_{k}\Vert _{\mathring{M}_{q,\lambda }\left( \Omega \right) }=1\quad
\text{for all }k  \label{eq seminorm 1}
\end{equation}%
and
\begin{equation}
\Vert \nabla f_{k}\Vert _{\mathring{M}_{q,\lambda }\left( \Omega \right)
}+\psi (f_{k})\rightarrow 0.  \label{eq seminorm}
\end{equation}%
By Theorem \ref{Theorem:RellichKondrachovZorkoSpaces}, we can suppose that $%
\Vert f_{k}-f\Vert _{\mathring{M}_{q,\lambda }\left( \Omega \right) }$ for
some $f\in \mathring{M}_{q,\lambda }\left( \Omega \right) .$ Moreover, by
\eqref{eq
seminorm}, we conclude that $\nabla f=0$ and $\psi (f)=0$. Then $f$ is
constant and $\psi (f)=|f|\psi (1)$ implies $f=0$, which contradicts
\eqref{eq
seminorm 1}. \fin

Our first Poincar\'{e}-type inequality takes the following form.

\begin{proposition}
\label{prop Poincare for f L1} Let $\Omega$ be a bounded and connected
Lipschitz domain in $\mathbb{R}^n$. Then, for all $f \in L^1(\Omega)$ with $%
\nabla f \in \mathring{M}_{q,\lambda}\left(\Omega\right)$, we have
\begin{equation*}
\|f\|_{\mathring{M}_{q,\lambda}\left(\Omega\right)} \leq c \left(\|\nabla
f\|_{\mathring{M}_{q,\lambda}\left(\Omega\right)}+\left| \int_\Omega f
\right|\right).
\end{equation*}
In particular, $f \in W^1\mathring{M}_{q,\lambda}\left(\Omega\right)$.
\end{proposition}

\noindent {\textbf{Proof.}} The proof is an adaptation of the one for Lemma
II.6.1 in \cite{galdi2011introduction}. Initially, let us suppose that $%
\Omega $ is \emph{star-shaped}, that is, there is $x_{0}\in \Omega $ such
that, for all $x\in \Omega $, the line segment from $x_{0}$ to $x$ is
contained in $\Omega $. By means of a translation, we can suppose $x_{0}=0$.
Now, let $\{r_{k}\}$ be an increasing sequence of positive numbers
converging to $1$ and, for each $k\in \mathbb{N}$, let
\begin{equation*}
\Omega _{k}:=\{x\in \mathbb{R}^{n}:r_{k}x\in \Omega \}.
\end{equation*}%
Then, $\Omega _{k}\supset \overline{\Omega }$. Let
\begin{equation*}
f_{k}(x):=f(r_{k}x),\quad x\in \Omega _{k},
\end{equation*}%
and, for $\epsilon (k)>0$ arbitrarily small, let $\phi _{k}:=\phi _{\epsilon
(k)}$ be a mollifier as in the Proposition \ref{Proposition:Mollifications}.
By properties of mollifiers in Lebesgue spaces, taken $\epsilon (k)<\mathrm{%
dist}(\Omega ,\mathbb{R}^{n}\backslash \Omega _{k})$, we can suppose that $%
\phi _{k}\ast f_{k}\in C^{\infty }(\overline{\Omega })$. In particular, it
follows from Corollary \ref{corollary before poincare} with $\psi
=|\int_{\Omega }\cdot \,|$ that, for all $j,k\in \mathbb{N}$,
\begin{equation}
\Vert \phi _{k}\ast f_{k}-\phi _{j}\ast f_{j}\Vert _{\mathring{M}_{q,\lambda
}\left( \Omega \right) }\leq c\left( \Vert \nabla (\phi _{k}\ast f_{k}-\phi
_{j}\ast f_{j})\Vert _{\mathring{M}_{q,\lambda }\left( \Omega \right)
}+\left\vert \int_{\Omega }(\phi _{k}\ast f_{k}-\phi _{j}\ast
f_{j})\right\vert \right) .  \label{ineq for phi ast f}
\end{equation}%
Since $\Vert \phi _{k}\Vert _{L^{1}(\mathbb{R}^{n})}=1$, by Young inequality
for convolutions we have
\begin{equation*}
\Vert \phi _{k}\ast f_{k}-f\Vert _{L^{1}(\Omega )}\leq \Vert \phi _{k}\ast
(f_{k}-f)\Vert _{L^{1}(\Omega )}+\Vert \phi _{k}\ast f-f\Vert _{L^{1}(\Omega
)}\leq c\Vert f_{k}-f\Vert _{L^{1}(\Omega )}+\Vert \phi _{k}\ast f-f\Vert
_{L^{1}(\Omega )}.
\end{equation*}%
Then $\phi _{k}\ast f_{k}\rightarrow f$ in $L^{1}(\Omega )$. Similarly, the
function
\begin{equation*}
\nabla (\phi _{k}\ast f_{k})=\phi _{k}\ast \nabla (f_{k})=(r_{k})^{n}\phi
_{k}\ast (\nabla f)_{k}
\end{equation*}%
converges to $\nabla f$ in $L^{n/\alpha }(\Omega )\subset \mathring{M}%
_{q,\lambda }\left( \Omega \right) $, where $\alpha :=(n-\lambda )/q$. Then, %
\eqref{ineq for phi ast f} shows us that $\{\phi _{k}\ast f_{k}\}$ is a
Cauchy sequence in $\mathring{M}_{q,\lambda }\left( \Omega \right) $.
Therefore, $f\in \mathring{M}_{q,\lambda }\left( \Omega \right) $ and
\begin{equation*}
\Vert f\Vert _{\mathring{M}_{q,\lambda }\left( \Omega \right) }\leq c\left(
\Vert \nabla f\Vert _{\mathring{M}_{q,\lambda }\left( \Omega \right)
}+\left\vert \int_{\Omega }f\right\vert \right) .
\end{equation*}

The general case in which $\Omega $ is just a bounded and connected
Lipschitz domain in $\mathbb{R}^{n}$ can be reduced to the previous one
since such domains are finite unions of star-shaped Lipschitz domains (see,
e.g., \cite[Exercise II.1.5]{galdi2011introduction}). \fin

Armed with the previous proposition, we can now obtain a version of the
Poincar\'{e} inequality in Zorko spaces.

\begin{theorem}[Poincar\'{e} Inequality for Zorko Spaces]
\label{Theorem:PoincareInequalityZorkoSpaces} Let $\Omega$ be a bounded and
connected Lipschitz domain in $\mathbb{R}^n$ and $U \subset \Omega$ with
positive measure. Then, there is $c>0$ such that
\begin{equation*}
\|f\|_{\mathring{M}_{q,\lambda}\left(\Omega\right)}\leq c\|\nabla f\|_{%
\mathring{M}_{q,\lambda}\left(\Omega\right)},
\end{equation*}
for all $f \in L^1(\Omega)$ satisfying $\nabla f \in \mathring{M}%
_{q,\lambda}\left(\Omega\right)$ and $\int_{U}f=0$, or $f \in W^1_0\mathring{%
M}_{q,\lambda}\left(\Omega\right)$.
\end{theorem}

\noindent {\textbf{Proof.}} For the case $\int_{U}f=0$, the inequality is a
direct consequence of Corollary \ref{corollary before poincare} with $\psi
(f):=\left\vert \int_{U}f\right\vert $ and Proposition \ref{prop Poincare
for f L1}. On the other hand, for $f\in W_{0}^{1}\mathring{M}_{q,\lambda
}\left( \Omega \right) $, the proof is as follows. First, let us suppose
that $f\in C_{0}^{\infty }(\Omega )$. Denoting $f_{0}=|\Omega
|^{-1}\int_{\Omega }f$, by the previous case, we have
\begin{equation}
\Vert f-f_{0}\Vert _{\mathring{M}_{q,\lambda }\left( \Omega \right) }\leq
c\Vert \nabla f\Vert _{\mathring{M}_{q,\lambda }\left( \Omega \right) }.
\label{inequality with mean}
\end{equation}%
Note that
\begin{equation*}
\Vert f_{0}\Vert _{\mathring{M}_{q,\lambda }\left( \Omega \right)
}=|f_{0}|\Vert 1\Vert _{\mathring{M}_{q,\lambda }\left( \Omega \right) }\leq
c\Vert f\Vert _{L^{q}(\Omega )}.
\end{equation*}%
By the classical Poincar\'{e} inequality in $W_{0}^{1,q}(\Omega )$, we have $%
\Vert f\Vert _{L^{q}(\Omega )}\leq c\Vert \nabla f\Vert _{L^{q}(\Omega )}$.
Furthermore, by Proposition \ref{embeddings from Kato} and since $|\Omega
|<\infty $, we have $\mathring{M}_{q,\lambda }\left( \Omega \right)
\hookrightarrow L^{q}(\Omega )$. Then,
\begin{equation*}
\Vert f_{0}\Vert _{\mathring{M}_{q,\lambda }\left( \Omega \right) }\leq
c\Vert \nabla f\Vert _{L^{q}(\Omega )}\leq c\Vert \nabla f\Vert _{\mathring{M%
}_{q,\lambda }\left( \Omega \right) },
\end{equation*}%
which, together with \eqref{inequality with mean}, implies $\Vert f\Vert _{%
\mathring{M}_{q,\lambda }\left( \Omega \right) }\leq c\Vert \nabla f\Vert _{%
\mathring{M}_{q,\lambda }\left( \Omega \right) }$. The general case follows
by density. \fin

The upcoming proposition pertains to a Bogovskii-type result within the
Zorko setting.

\begin{proposition}[Bogovskii-type Proposition for Zorko Spaces]
\label{proposition Bogovskii} Let $\Omega \subset \mathbb{R}^{n}$ be a
bounded Lipschitz domain. Then, given $f\in \mathring{M}_{q,\lambda }\left(
\Omega \right) $, satisfying $\int_{\Omega }f=0,$ there exists a solution $%
\mathbf{w}\in W_{0}^{1}\mathring{M}_{q,\lambda }\left( \Omega \right) $ for
\begin{equation}
\nabla \cdot \mathbf{w}(x)=f(x),\quad x\in \Omega ,  \label{eq div}
\end{equation}%
such that
\begin{equation}
\Vert \mathbf{w}\Vert _{W_{0}^{1}\mathring{M}_{q,\lambda }\left( \Omega
\right) }\leq c\Vert f\Vert _{\mathring{M}_{q,\lambda }\left( \Omega \right)
},  \label{estimate div}
\end{equation}%
where $c=c(n,q,\lambda ,\Omega )$ is independent of $f$.
\end{proposition}

\noindent {\textbf{Proof.}} The proof is closely based on the arguments of
\cite[Section III.3]{galdi2011introduction}. We restrict ourselves to those
ones that need to be adapted.

First, we consider a especial type of domain, which we will call \emph{%
ball-star-shaped domain}, meaning that there is a ball $\overline{B(x_{0},R)}%
\subset \Omega $ such that $\Omega $ is star-shaped with respect to every
point in $B(x_{0},R)$, that is, for all $x\in B(x_{0},R)$ and $y\in \Omega $%
, the line segment from $x$ to $y$ is contained in $\Omega $. By a suitable
change of variables, we can suppose $R=1$ and $x_{0}=0$. We also consider
initially that $f\in C_{0}^{\infty }(\Omega )$. In \cite%
{galdi2011introduction}, Galdi proved that if $\omega \in C_{0}^{\infty
}(\Omega )$ is chosen with $\mathrm{supp}(\omega )\subset B_{1}$ and $%
\int_{\Omega }\omega =1$, and we define
\begin{equation*}
W(x,y):=(x-y)\int_{1}^{\infty }\omega (y+r(x-y))r^{n-1}dr,
\end{equation*}%
then the vector field
\begin{equation*}
\mathbf{w}(x)=\int_{\Omega }W(x,y)f(y)dy,
\end{equation*}%
is smooth with compact support in $\Omega $ and solves \eqref{eq div}.
Moreover, for $j,i=1,2,...,n$,
\begin{gather*}
\partial _{j}v_{i}(x)=\int_{\Omega }\frac{\nu _{ij}(x,x-y)}{|x-y|^{n}}%
f(y)dy+\int_{\Omega }G_{ij}(x,y)f(y)dy+f(x)\int_{\Omega }\frac{%
(x_{j}-y_{j})(x_{i}-y_{i})}{|x-y|^{2}}dy \\
=f_{1}(x)+f_{2}(x)+f_{3}(x),
\end{gather*}%
where $\nu _{ij}$ and $G_{ij}$ are certain functions satisfying
\begin{equation*}
\nu _{ij}(x,y)=\nu _{ij}(x,\alpha y),\quad \alpha >0,\,x\in \Omega ,y\in
\mathbb{R}^{n}\backslash {0},
\end{equation*}%
\begin{equation*}
\int_{|y|=1}\nu _{ij}(x,y)dy=0,
\end{equation*}%
\begin{equation*}
|\nu _{ij}(x,y)|\leq C,\quad x\in \Omega ,\,|y|=1,
\end{equation*}%
and
\begin{equation*}
|G_{ij}(x,y)|\leq c|x-y|^{1-n},\quad x,y\in \Omega .
\end{equation*}%
By Proposition \ref{Proposition Calderon-Zygmund} (Calder\'{o}n-Zygmund
inequality in $M_{q,\lambda }\left( \Omega \right) $),
\begin{equation*}
\Vert f_{1}\Vert _{M_{q,\lambda }\left( \Omega \right) }\leq c\Vert f\Vert _{%
\mathring{M}_{q,\lambda }\left( \Omega \right) },
\end{equation*}%
and from Proposition \ref{proposition fractional integral operators} it
follows that
\begin{equation*}
\Vert f_{2}\Vert _{M_{q,\lambda }\left( \Omega \right) }\leq c\Vert f\Vert _{%
\mathring{M}_{q,\lambda }\left( \Omega \right) }.
\end{equation*}%
Of course,
\begin{equation*}
\Vert f_{3}\Vert _{M_{q,\lambda }\left( \Omega \right) }\leq c\Vert f\Vert _{%
\mathring{M}_{q,\lambda }\left( \Omega \right) }.
\end{equation*}%
hence
\begin{equation*}
\Vert \nabla \mathbf{w}\Vert _{\mathring{M}_{q,\lambda }\left( \Omega
\right) }\leq c\Vert f\Vert _{\mathring{M}_{q,\lambda }\left( \Omega \right)
}
\end{equation*}%
Therefore, by the Poincar\'{e} inequality for $W_{0}^{1}\mathring{M}%
_{q,\lambda }\left( \Omega \right) $ (Theorem \ref%
{Theorem:PoincareInequalityZorkoSpaces}), we have
\begin{equation*}
\Vert \mathbf{w}\Vert _{W_{0}^{1}\mathring{M}_{q,\lambda }\left( \Omega
\right) }\leq c\Vert f\Vert _{\mathring{M}_{q,\lambda }\left( \Omega \right)
}.
\end{equation*}%
Since the map $f\mapsto \mathbf{w}$ constructed above is linear, the result
remains valid for any $f\in \mathring{M}_{q,\lambda }\left( \Omega \right) $
by an argument of density.

For the general case, we use the fact that a bounded Lipschitz domain $%
\Omega $ is a finite union $\cup_1^N \Omega_k$, where each $\Omega_k$ is a
ball-star-shaped domain (see \cite{galdi2011introduction}, Lemma II.1.3).
Then, given $f \in \mathring{M}_{q,\lambda}\left(\Omega\right)$ such that $%
\int_\Omega f=0$, we construct functions $f_k$, $k=1,2,...,N$ satisfying

\begin{itemize}
\item[(I)] $\mathrm{supp}(f_k) \subset \overline{\Omega_k}$,

\item[(II)] $\int_{\Omega_k}f_k(x)\,dx=0,$

\item[(III)] $f=\sum_1^N f_k$.
\end{itemize}

Then, by solving $\nabla \cdot \mathbf{w}_{k}=f_{k}$ in each $\Omega _{k}$
and setting $\mathbf{w}=\sum_{1}^{N}\mathbf{w}_{k}$ in $\Omega $, it is
sufficient to prove that
\begin{equation}
\sum_{1}^{N}\Vert f_{k}\Vert _{\mathring{M}_{q,\lambda }\left( \Omega
\right) }\leq c\Vert f\Vert _{\mathring{M}_{q,\lambda }\left( \Omega \right)
}.  \label{inequality for lip domain in bogovskii}
\end{equation}%
The functions are constructed as follows. For $k=1,2,...,N-1$, let
\begin{equation*}
D_{k}=\cup _{j=k+1}^{N}\Omega _{j},\quad F_{k}=\Omega _{k}\cap D_{k},
\end{equation*}%
and let $\chi _{k}$ be the characteristic function of $F_{k}$. Then, we set $%
g_{0}=f$ and, for $k=1,2,...,N-1$,
\begin{equation*}
g_{k}(x)=(1-\chi _{k}(x))g_{k-1}(x)-\frac{\chi _{k}(x)}{|F_{k}|}%
\int_{D_{k}\backslash \Omega _{k}}g_{k-1},\quad x\in D_{k}.
\end{equation*}%
Finally, for $k=1,2,...,N-1$,
\begin{equation*}
f_{k}(x)=g_{k-1}(x)-\frac{\chi _{k}(x)}{|F_{k}|}\int_{\Omega
_{k}}g_{k-1},\quad x\in \Omega _{k},
\end{equation*}%
and $f_{N}=g_{N-1}$. It is not hard to verify (I), (II) and (III). Moreover,
for $k=1,2,...N-1$, it holds estimates of the type
\begin{equation*}
\Vert f_{k}\Vert _{\mathring{M}_{q,\lambda }\left( \Omega \right) }\leq
c_{k}\Vert g_{k-1}\Vert _{\mathring{M}_{q,\lambda }\left( \Omega \right)
},\quad k=1,2,...,N,
\end{equation*}%
and
\begin{equation*}
\Vert g_{k}\Vert _{\mathring{M}_{q,\lambda }\left( \Omega \right) }\leq
c_{k}^{\prime }\Vert g_{k-1}\Vert _{\mathring{M}_{q,\lambda }\left( \Omega
\right) },\quad k=1,2,...,N-1,
\end{equation*}%
from what we conclude that $\Vert f_{k}\Vert _{\mathring{M}_{q,\lambda
}\left( \Omega \right) }\leq c\Vert g_{0}\Vert _{\mathring{M}_{q,\lambda
}\left( \Omega \right) }=c\Vert f\Vert _{\mathring{M}_{q,\lambda }\left(
\Omega \right) }$, hence \eqref{inequality for
lip domain in bogovskii} holds. \fin

We conclude this section by presenting the Rellich-Kondrachov theorem and
the Poincar\'{e} inequality in the context of block spaces. The content of
the next theorem addresses the first of these results.

\begin{theorem}[Rellich-Kondrachov Theorem for Block Spaces]
Let $\Omega \subset \mathbb{R}^n$ be a bounded Lipschitz domain. Then the
following compact embedding holds:
\begin{equation}  \label{compact embedding block}
W^1 H_{q,\lambda}(\Omega) \hookrightarrow\hookrightarrow
H_{q,\lambda}(\Omega).
\end{equation}
\end{theorem}

\noindent{\textbf{Proof.}} Let $\{f_k\}$ be a bounded sequence in $%
W^1H_{q,\lambda}(\Omega)$. By Banach Alaoglu theorem, we can suppose that
there is $f \in W^1H_{q,\lambda}(\Omega)$ such that $f_k
\rightharpoonup^\ast f$ and $\nabla h_k \rightharpoonup^\ast \nabla f$ in
the weak-star topology of $H_{q,\lambda}(\Omega)=\mathring{M}%
_{q^{\prime},\lambda}(\Omega)^\ast$. We claim that $f_k \to f $ in ${%
H_{q,\lambda}(\Omega)}$. Indeed, let us suppose by contradiction that, up to
a subsequence, $\|f_k - f\|_{H_{q,\lambda}(\Omega)}\geq 2\epsilon$ for some $%
\epsilon>0$. Then, for each $k$, there is $g_k \in \mathring{M}%
_{q^{\prime},\lambda}\left(\Omega\right)$ such that $\|g_k\|_{\mathring{M}%
_{q^{\prime},\lambda}\left(\Omega\right)}=1$ and
\begin{equation*}
\int_\Omega (f_k-f)g_k \geq \epsilon.
\end{equation*}
Denoting $h_k:=g_k - d_k$, where $d_k:=|\Omega|^{-1}\int_\Omega g_k $, by
Proposition \ref{proposition Bogovskii}, there is a vector field $\mathbf{v}%
_k \in W^1_0\mathring{M}_{q^{\prime},\lambda}\left(\Omega\right)$ such that $%
\|\mathbf{v}_k\|_{W^1_0\mathring{M}_{q^{\prime},\lambda}\left(\Omega\right)}
\leq c \|h_k\|_{\mathring{M}_{q^{\prime},\lambda}\left(\Omega\right)}$ and
\begin{equation*}
\nabla \cdot \mathbf{v}_k=h_k.
\end{equation*}
Note that $\|h_k\|_{\mathring{M}_{q^{\prime},\lambda}\left(\Omega\right)}%
\leq c \|g_k\|_{\mathring{M}_{q^{\prime},\lambda}\left(\Omega\right)} \leq c$%
. Then $\{\mathbf{v}_k\}$ is a bounded sequence in $W^1_0 \mathring{M}%
_{q^{\prime},\lambda}\left(\Omega\right)$ and by Rellich-Kondrachov theorem
for Zorko spaces there is $\mathbf{v} \in W^1 \mathring{M}%
_{q^{\prime},\lambda}\left(\Omega\right)$ such that, up to a subsequence, $\|%
\mathbf{v}_k - \mathbf{v}\|_{\mathring{M}_{q^{\prime},\lambda}\left(\Omega%
\right)}\to 0.$ Therefore,
\begin{align*}
\epsilon \leq \int_\Omega (f_k-f)g_k = \int_\Omega (f_k-f)(\nabla\cdot
\mathbf{v}_k + d_k) = \int_\Omega \nabla(f_k - f)\cdot \mathbf{v}_k +
d_k\int_\Omega (f_k - f),
\end{align*}
which is a contradiction, since $\nabla f_k \rightharpoonup^\ast \nabla f$, $%
\mathbf{v}_k \to \mathbf{v}$, $|d_k|\leq c\|g_k\|_{\mathring{M}%
_{q^{\prime},\lambda}\left(\Omega\right)}\leq c$ and $f_k
\rightharpoonup^\ast f$. \fin

The upcoming theorem presents the statements of the analogous versions of
Corollary \ref{corollary before poincare}, Proposition \ref{prop Poincare
for f L1}, and Theorem \ref{Theorem:PoincareInequalityZorkoSpaces} for block
spaces. Its proof closely resembles the previous ones and will therefore be
omitted.

\begin{theorem}[Poincar\'{e} Inequality for Block Spaces]
\label{theorem Poincare inequality block} Let $\Omega$ be a bounded and
connected Lipschitz domain in $\mathbb{R}^n$. Then the following
propositions hold.

\begin{itemize}
\item[(I)] Let $\psi:H_{q,\lambda}\left(\Omega\right)\to \mathbb{R}_+$ be a
continuous and absolutely homogeneous function ($\psi(sf)=|s|\psi(f)$) such
that $\psi(1)\neq 0$. Then there is a constant $c$ such that
\begin{equation*}
\|f\|_{\mathring{M}_{q,\lambda}\left(\Omega\right)}\leq c \left(\|\nabla f
\|_{\mathring{M}_{q,\lambda}\left(\Omega\right)}+\psi(f)\right),
\end{equation*}
for all $f \in W^1H_{q,\lambda}\left(\Omega\right)$.

\item[(II)] There is $c>0$ such that, for all $f \in L^1(\Omega)$ with $%
\nabla f \in H_{q,\lambda}\left(\Omega\right)$, we have
\begin{equation*}
\|f\|_{\mathring{M}_{q,\lambda}\left(\Omega\right)} \leq c \left(\|\nabla
f\|_{\mathring{M}_{q,\lambda}\left(\Omega\right)}+\left| \int_\Omega f
\right|\right).
\end{equation*}
In particular, $f \in W^1H_{q,\lambda}\left(\Omega\right)$.

\item[(III)] There is $c>0$ such that
\begin{equation*}
\|f\|_{H_{q,\lambda}\left(\Omega\right)}\leq c\|\nabla
f\|_{H_{q,\lambda}\left(\Omega\right)},
\end{equation*}
for all $f \in L^1(\Omega)$ satisfying $\nabla f \in
H_{q,\lambda}\left(\Omega\right)$ and $\int_{U}f=0$, or $f \in
W^1_0H_{q,\lambda}\left(\Omega\right)$.
\end{itemize}
\end{theorem}

\section{Proof of the Main Result (Theorem \protect\ref%
{theorem:HelmholtzDecomposition})}

\label{section main result}

This section is dedicated to proving our main result, structured into three
subsections related to the types of domain $\Omega $ and function spaces.

We begin by revisiting a de Rham lemma, the proof of which can be found in
\cite[Lemma III.1.1.]{galdi2011introduction}.

\begin{lemma}[De Rham]
\label{lema de Rham} Let $\Omega$ be an non-empty open subset of $\mathbb{R}%
^n$ and $\mathbf{v}$ a vector field in $L^1_{\mathrm{loc}}(\Omega)$. Suppose
that, for all $\mathbf{w}\in C^\infty_{0,\sigma}(\Omega)$, we have $%
\int_\Omega \mathbf{v}\cdot \mathbf{w}=0.$ Then, $\mathbf{v}=\nabla p $ for
some $p \in W^{1,1}_{\mathrm{loc}}(\overline{\Omega})$.
\end{lemma}

\begin{remark}
\label{remark irrotational vector fields are closed} As a direct consequence
of Lemma \ref{lema de Rham}, we obtain that the spaces of irrotational
vector fields are complete. Indeed, given a sequence $\{\nabla p_k\} \subset
GM_{q,\lambda}\left(\Omega\right)$ and $\mathbf{v} \in
M_{q,\lambda}\left(\Omega\right)$ such that $\|\nabla p_k-\mathbf{v}%
\|_{M_{q,\lambda}\left(\Omega\right)}\to 0$, we have
\begin{equation*}
\int_{\Omega} \mathbf{v}\cdot \mathbf{w}=\lim \int_{\Omega} \nabla p_k\cdot
\mathbf{w} = -\lim\int_{\Omega} p_k\nabla\cdot \mathbf{w} = 0,
\end{equation*}
for all for all $\mathbf{w}\in C^\infty_{0,\sigma}(\Omega)$. Then $v=\nabla
p $ for some $p \in W^{1,1}_{\mathrm{loc}}(\overline{\Omega})$, that is, $%
\mathbf{v} \in GM_{q,\lambda}\left(\Omega\right)$. Of course the same
argument holds for $G\mathring{M}_{q,\lambda}\left(\Omega\right)$ and $%
GH_{q,\lambda}(\Omega)$.
\end{remark}

The Helmholtz decomposition $\nabla p+(\mathbf{u}-\nabla p)$ in $%
GL^{q}(\Omega )\oplus SL^{q}(\Omega )$ of a vector field $\mathbf{u}\in
L^{q}(\Omega )$ is equivalent to the well-posedness of a weak version of the
Neumann problem (see \cite{simader1992new})
\begin{align*}
\Delta p& =\nabla \cdot \mathbf{u},\quad \text{in }\Omega , \\
\frac{\partial p}{\partial \mathbf{n}}& =\mathbf{u}\cdot \mathbf{n},\quad
\text{on }\partial \Omega .
\end{align*}%
The next lemma generalizes this result to the setting of Zorko subspaces.

\begin{lemma}[Equivalence Lemma]
\label{lemma de equivalencia} Let $\Omega$ be a non-empty open subset of $%
\mathbb{R}^n$, $\mathbf{u} \in \mathring{M}_{q,\lambda}\left(\Omega\right)$
and $\nabla p \in G\mathring{M}_{q,\lambda}\left(\Omega\right)$. Then, $%
\mathbf{u}=\nabla p + (\mathbf{u}-\nabla p)$ is a Helmholtz decomposition of
$\mathbf{u}$ in $\mathring{M}_{q,\lambda}\left(\Omega\right)$, that is, $%
\mathbf{u}-\nabla p$ belongs to $S\mathring{M}_{q,\lambda}\left(\Omega%
\right) $, if and only if $p$ satisfies
\begin{equation}  \label{weak neumann problem}
\int_{\Omega}\nabla p\cdot\nabla \phi=\int_{\Omega}\mathbf{u}\cdot\nabla
\phi, \quad \text{for all }\nabla \phi\in
GH_{q^{\prime},\lambda}\left(\Omega\right).
\end{equation}
In particular, the Helmholtz decomposition problem in $\mathring{M}%
_{q,\lambda}\left(\Omega\right)$ is equivalent to prove that, for each $%
\mathbf{u} \in \mathring{M}_{q,\lambda}\left(\Omega\right)$, there is a
unique solution $\nabla p \in G\mathring{M}_{q,\lambda}\left(\Omega\right)$
to \eqref{weak neumann problem}, and such solution satisfies
\begin{equation*}
\|\nabla p\|_{\mathring{M}_{q,\lambda}\left(\Omega\right)}\leq c\|\mathbf{u}%
\|_{\mathring{M}_{q,\lambda}\left(\Omega\right)},
\end{equation*}
for some $c>0$ independent of $\mathbf{u}$.
\end{lemma}

\noindent{\textbf{Proof.}} Let us denote
\begin{equation*}
\mathbf{w}=\mathbf{u}-\nabla p \in \mathring{M}_{q,\lambda}\left(\Omega%
\right).
\end{equation*}
If $\mathbf{w} \in S\mathring{M}_{q,\lambda}\left(\Omega\right)$, then $%
\mathbf{w}$ can be approached in $\mathring{M}_{q,\lambda}\left(\Omega%
\right) $ by functions in $C^\infty_{0,\sigma}(\Omega)$, so
\eqref{weak
neumann problem} holds. Reciprocally, if $\mathbf{w} \not\in S\mathring{M}%
_{q,\lambda}\left(\Omega\right)$, by Hahn-Banach theorem, there is $\mathbf{v%
} \in H_{q^{\prime},\lambda}\left(\Omega\right)$ such that $\int_{\Omega}%
\mathbf{w}\cdot \mathbf{v}\neq 0$ but
\begin{equation}  \label{eq Hanh Banach}
\int_{\Omega}\mathbf{w}_0\cdot \mathbf{v}= 0,\quad \text{for all } \mathbf{w}%
_0 \in S\mathring{M}_{q,\lambda}\left(\Omega\right).
\end{equation}
By Lemma \ref{lema de Rham}, \eqref{eq Hanh Banach} implies that $\mathbf{v}
\in GH_{q^{\prime},\lambda}\left(\Omega\right)$. Then
\eqref{weak neumann
problem} does not hold. \fin

\subsection{Helmholtz Decomposition for $\mathring{M}_{q,\protect\lambda%
}\left(\mathbb{R}^n\right)$ and $\mathring{M}_{q,\protect\lambda}\left(%
\mathbb{R}^n_+\right)$}

\label{subsection Helmholtz Decomposition for Zorko spaces in half space}

Here, we address the simplest cases of $\mathbb{R}^{n}$ and $\mathbb{R}%
_{+}^{n}.$ We begin with $\mathring{M}_{q,\lambda }\left( \mathbb{R}%
^{n}\right) $. Given $\mathbf{u}\in C_{0}^{\infty }(\mathbb{R}^{n})^{n}$,
the Helmholtz decomposition in $L^{\frac{n}{\alpha }}(\mathbb{R}^{n})$,
where $\alpha =(n-\lambda )/q$, asserts that there exist unique vector
fields
\begin{equation*}
\mathbf{v}\in GL^{\frac{n}{\alpha }}(\mathbb{R}^{n})\subset G\mathring{M}%
_{q,\lambda }\left( \mathbb{R}^{n}\right) ,\quad \mathbf{w}\in SL^{\frac{n}{%
\alpha }}(\mathbb{R}^{n})\subset S\mathring{M}_{q,\lambda }\left( \mathbb{R}%
^{n}\right) ,
\end{equation*}%
such that $\mathbf{u}=\mathbf{v}+\mathbf{w}$. Moreover, using the Calder\'{o}%
n-Zygmund theorem and Lemma \ref{lemma de equivalencia} for Lebesgue spaces,
it is not difficult to verify that $\mathbf{v}=\nabla p$ with $p$ given by
\begin{equation*}
p(x):=\int_{\mathbb{R}^{n}}\Gamma (x-y)\nabla \cdot \mathbf{u}(y)\,dy=\int_{%
\mathbb{R}^{n}}(\nabla \Gamma )(x-y)\cdot \mathbf{u}(y)\,dy,\quad x\in
\mathbb{R}^{n},
\end{equation*}%
(see \cite{galdi2011introduction}, Section III.1). We remind that $\Gamma $
denotes the fundamental solution for the Laplacian in $\mathbb{R}^{n}$. For $%
\epsilon >0$, let
\begin{equation*}
p_{\epsilon }(x):=\int_{|x-y|>\epsilon }(\nabla \Gamma )(x-y)\cdot \mathbf{u}%
(y)\,dy.
\end{equation*}%
Then, $p_{\epsilon }\rightarrow p$ pointwise as $\epsilon \rightarrow 0$.
Moreover,
\begin{equation*}
\nabla p_{\epsilon }(x)=\int_{|x-y|>\epsilon }(\nabla ^{2}\Gamma )(x-y)\cdot
\mathbf{u}(y)\,dy\quad +\int_{|x-y|=\epsilon }(\nabla \Gamma )(x-y)\cdot
\mathbf{u}(y)\frac{x-y}{\epsilon }\,d\sigma _{y},
\end{equation*}%
where $\nabla ^{2}\Gamma $ is the Hessian matrix of $\Gamma $ and $\sigma
_{y}$ is the standard surface measure on the sphere $\{y\in \mathbb{R}%
^{n}:|x-y|=\epsilon \}$. Since $\mathbf{u}\in C_{0}^{\infty }(\mathbb{R}%
^{n}) $, we obtain
\begin{equation*}
\lim_{\epsilon \rightarrow 0}\nabla p_{\epsilon }(x)=\text{p.v.}\int_{%
\mathbb{R}^{n}}(\nabla ^{2}\Gamma )(x-y)\cdot \mathbf{u}(y)\,dy+\frac{%
\mathbf{u}(x)}{n},
\end{equation*}%
with the convergence of the limit being uniform for $x$ in compact sets, and
the integral is understood in the sense of the principal value. Therefore,
\begin{equation*}
\nabla p(x)=\text{p.v.}\int_{\mathbb{R}^{n}}(\nabla ^{2}\Gamma )(x-y)\cdot
\mathbf{u}(y)\,dy+\frac{\mathbf{u}(x)}{n},
\end{equation*}%
which leads to the estimate
\begin{equation}
\Vert \nabla p\Vert _{\mathring{M}_{q,\lambda }\left( \Omega \right) }\leq
c\Vert \mathbf{u}\Vert _{\mathring{M}_{q,\lambda }\left( \Omega \right) },
\label{estimate nabla p in Rn}
\end{equation}%
by Proposition \ref{Proposition Calderon-Zygmund}. The general case $\mathbf{%
u}\in \mathring{M}_{q,\lambda }\left( \mathbb{R}^{n}\right) $ follows by
density, since $G\mathring{M}_{q,\lambda }\left( \mathbb{R}^{n}\right) $ is
closed.

For $\mathring{M}_{q,\lambda}\left(\mathbb{R}^n_+\right)$, the proof
proceeds in a similar manner. In this case, for $\mathbf{u} \in C^\infty_0(%
\mathbb{R}^n_+)^n$, we have
\begin{align*}
p(x) & := -\int_{\mathbb{R}^n_+} N(x, y) \nabla \cdot \mathbf{u}(y) \, dy,
\end{align*}
where $N$ is the Green (Neumann) function for the Laplacian in $\mathbb{R}%
^n_+$, that is,
\begin{equation*}
N(x, y) = \Gamma(x - \overline{y}) - \Gamma(x - y), \quad x, y \in \mathbb{R}%
^n_+,
\end{equation*}
with $\overline{y} = (y_1, \ldots, y_{n-1}, -y_n)$.

\subsection{Helmholtz Decomposition for $\mathring{M}_{q,\protect\lambda%
}\left(\Omega\right)$ with Bounded or Exterior $C^1$ Domains}

\label{subsection Helmholtz Decomposition for Zorko spaces in bounded or
exterior domains}

In this section, we establish the Helmholtz decomposition for $\mathring{M}%
_{q,\lambda}(\Omega)$, where $\Omega$ is either a bounded or an exterior $%
C^1 $ domain. To achieve this, we utilize the equivalent formulation
provided by the weak Neumann problem in Lemma \ref{lemma de equivalencia}.
The following lemma serves as a uniqueness result.

\begin{lemma}
\label{coerciveness lemma} Let $\Omega \subset \mathbb{R}^{n}$ be either a
bounded or exterior $C^{1}$ domain. If $\nabla p\in G\mathring{M}_{q,\lambda
}\left( \Omega \right) $ satisfies
\begin{equation*}
\int_{\Omega }\nabla p\cdot \nabla \phi =0
\end{equation*}%
for all $\nabla \phi \in GH_{q^{\prime },\lambda }\left( \Omega \right) $,
then $\nabla p=0$.
\end{lemma}

\noindent {\textbf{Proof.}} The claim holds for weighted Lebesgue spaces $%
L_{w}^{q}(\Omega )$ with a weight $w$ in the Muckenhoupt class $A_{q}$ (see
Lemma 2 in \cite{frohlich2000helmholtz}), that is, with $\mathring{M}%
_{q,\lambda }\left( \Omega \right) $ and $H_{q^{\prime },\lambda }\left(
\Omega \right) $ replaced by $L_{w}^{q}(\Omega )$ and $L_{w^{\prime
}}^{q^{\prime }}(\Omega )$, respectively. Then the proof follows from
Proposition \ref{embeddings from Kato}, which states that $\mathring{M}%
_{q,\lambda }\left( \Omega \right) \hookrightarrow \L _{w}^{q}(\Omega )$ and
$H_{q^{\prime },\lambda }\left( \Omega \right) \hookleftarrow L_{w^{\prime
}}^{q^{\prime }}(\Omega )$ for a suitable weight $w\in A_{q}$. \fin

Any vector field $\mathbf{u}\in \mathring{M}_{q,\lambda }\left( \Omega
\right) $ can be interpreted as a continuous linear functional in $%
GH_{q^{\prime },\lambda }\left( \Omega \right) ^{\ast }$ defined by the
mapping $\nabla \phi \mapsto \langle \mathbf{u},\nabla \phi \rangle
:=\int_{\Omega }\mathbf{u}\cdot \nabla \phi $ for $\nabla \phi \in
GH_{q^{\prime },\lambda }\left( \Omega \right) $. Then, the weak Neumann
problem described in Lemma \ref{lemma de equivalencia} consists of proving
that any such $\mathbf{u}$, regarded as an element of $GH_{q^{\prime
},\lambda }\left( \Omega \right) ^{\ast }$, can be uniquely represented by a
gradient field $\nabla p\in G\mathring{M}_{q,\lambda }\left( \Omega \right) $%
. We will solve this by showing that the map $\nabla p\mapsto \langle \nabla
p,\cdot \rangle $ is indeed an isomorphism from $G\mathring{M}_{q,\lambda
}\left( \Omega \right) $ to $GH_{q^{\prime },\lambda }\left( \Omega \right)
^{\ast }$. The first step is to demonstrate that it has a closed range,
which follows from the variational inequality stated in the theorem below.
Additionally, this inequality may serve as a tool for obtaining estimates of
gradient fields in the Morrey setting, which could be of independent
interest.

\begin{theorem}[Variational Inequality]
\label{variational inequality theorem} Let $\Omega $ be a bounded or
exterior $C^{1}$ domain. Then there is $c>0$ such that
\begin{equation}
\Vert \nabla p\Vert _{\mathring{M}_{q,\lambda }\left( \Omega \right) }\leq
c\sup \left\{ \frac{\int_{\Omega }\nabla p\cdot \nabla \phi }{\Vert \nabla
\phi \Vert _{H_{q^{\prime },\lambda }\left( \Omega \right) }}\right\} ,
\label{variational ineq in statement}
\end{equation}%
for all $\nabla p\in G\mathring{M}_{q,\lambda }\left( \Omega \right) $,
where the supremum is taken over all non-null functions $\nabla \phi \in
GH_{q^{\prime },\lambda }\left( \Omega \right) $.
\end{theorem}

\noindent {\textbf{Proof.}} We start by observing that, by duality and the
Helmholtz decomposition for $H_{q^{\prime },\lambda }\left( \mathbb{R}%
_{+}^{n}\right) $, it is easy to verify that the inequality holds if $\Omega
$ is replaced by $\mathbb{R}_{+}^{n}$. This allows us to extend the result
to a $C^{1}$ slightly perturbed half-space in the following sense. Given $%
\sigma :\mathbb{R}^{n-1}\rightarrow \mathbb{R}$ in $C_{0}^{1}(\mathbb{R}%
^{n-1})$ such that $|\nabla \sigma (0)|=\sigma (0)=0$, let
\begin{equation*}
S_{\sigma }:=\{x=(x^{\prime },x_{n})\in \mathbb{R}^{n}:x_{n}>\sigma
(x^{\prime })\}.
\end{equation*}%
Then the change of coordinates $y(x)=(x^{\prime },x_{n}-\sigma (x^{\prime
})) $ is a $C^{1}$ diffeomorphism between $S_{\sigma }$ and $\mathbb{R}%
_{+}^{n}$ with Jacobian matrix
\begin{equation*}
\frac{dy}{dx}=\left[
\begin{matrix}
I_{n-1} & 0 \\
-\nabla \sigma (x^{\prime }) & 1%
\end{matrix}%
\right] .
\end{equation*}%
Since $\det (dy/dx)=1$, the mapping $f\mapsto Kf$ defined by $Kf(y):=f(x(y))$
is an isomorphism from $\mathring{M}_{q,\lambda }\left( S_{\sigma }\right) $
to $\mathring{M}_{q,\lambda }\left( \mathbb{R}_{+}^{n}\right) $ and from $%
H_{q^{\prime },\lambda }\left( S_{\sigma }\right) $ to $H_{q^{\prime
},\lambda }\left( \mathbb{R}_{+}^{n}\right) $. Moreover, It is not hard to
verify that, for $\nabla p\in G\mathring{M}_{q,\lambda }\left( S_{\sigma
}\right) $ and $\nabla \phi \in GH_{q^{\prime },\lambda }\left( S_{\sigma
}\right) $, we have
\begin{equation*}
\Vert \nabla p\Vert _{\mathring{M}_{q,\lambda }\left( S_{\sigma }\right)
}\leq c(1+R)\Vert \nabla (Kp)\Vert _{\mathring{M}_{q,\lambda }\left( \mathbb{%
R}_{+}^{n}\right) },
\end{equation*}%
\begin{equation*}
\Vert \nabla \phi \Vert _{H_{q^{\prime },\lambda }\left( S_{\sigma }\right)
}\leq c(1+R)\Vert \nabla (K\phi )\Vert _{H_{q^{\prime },\lambda }\left(
\mathbb{R}_{+}^{n}\right) },
\end{equation*}%
and
\begin{equation*}
\left\vert \int_{\mathbb{R}_{+}^{n}}\nabla (Kp)\cdot \nabla (K\phi
)\,dy\right\vert \leq \left\vert \int_{S_{\sigma }}\nabla p\cdot \nabla \phi
\,dx\right\vert +cR(1+R)\Vert \nabla p\Vert _{\mathring{M}_{q,\lambda
}\left( S_{\sigma }\right) }\Vert \nabla \phi \Vert _{H_{q^{\prime },\lambda
}\left( S_{\sigma }\right) },
\end{equation*}%
where $R:=\Vert \nabla \sigma \Vert _{L^{\infty }(\mathbb{R}_{+}^{n-1})}$
and $c>0$ is some constant independent of $p,\phi $ or $\sigma $. Then,
\begin{equation*}
\frac{\left\vert \int_{\mathbb{R}_{+}^{n}}\nabla (Kp)\cdot \nabla (K\phi
)\,dy\right\vert }{\Vert \nabla (K\phi )\Vert _{H_{q^{\prime },\lambda
}\left( \mathbb{R}_{+}^{n}\right) }}\leq c(1+R)\frac{\left\vert
\int_{S_{\sigma }}\nabla p\cdot \nabla \phi \,dx\right\vert }{\Vert \nabla
\phi \Vert _{H_{q^{\prime },\lambda }\left( S_{\sigma }\right) }}%
+cR(1+R)^{2}\Vert \nabla p\Vert _{\mathring{M}_{q,\lambda }\left( S_{\sigma
}\right) }.
\end{equation*}%
Taking the supremum over $\nabla \phi \in GH_{q^{\prime },\lambda }\left(
S_{\sigma }\right) $, we have
\begin{equation*}
\Vert \nabla (Kp)\Vert _{\mathring{M}_{q,\lambda }\left( \mathbb{R}%
_{+}^{n}\right) }\leq c(1+R)\sup \left\{ \frac{\left\vert \int_{S_{\sigma
}}\nabla p\cdot \nabla \phi \,dx\right\vert }{\Vert \nabla \phi \Vert
_{H_{q^{\prime },\lambda }\left( S_{\sigma }\right) }}\right\}
+cR(1+R)^{2}\Vert \nabla p\Vert _{\mathring{M}_{q,\lambda }\left( S_{\sigma
}\right) }.
\end{equation*}%
Hence,
\begin{equation*}
\Vert \nabla p\Vert _{\mathring{M}_{q,\lambda }\left( S_{\sigma }\right)
}\leq c\sup \left\{ \frac{\left\vert \int_{S_{\sigma }}\nabla p\cdot \nabla
\phi \,dx\right\vert }{\Vert \nabla \phi \Vert _{H_{q^{\prime },\lambda
}\left( S_{\sigma }\right) }}\right\} +cR(1+R)\Vert \nabla p\Vert _{%
\mathring{M}_{q,\lambda }\left( S_{\sigma }\right) }.
\end{equation*}%
Therefore, if $R=\Vert \nabla \sigma \Vert _{L^{\infty }(\mathbb{R}%
_{+}^{n-1})}$ is small enough, then inequality
\eqref{variational ineq in
statement} holds with $\Omega =S_{\sigma }$.

Now, we prove \eqref{variational ineq in statement} in the case in which $%
\Omega$ is a bounded or exterior domain using an argument of contradiction.
Suppose by contradiction that there is a sequence $\{\nabla p_k\}\subset G%
\mathring{M}_{q,\lambda}\left(\Omega\right)$ such that
\begin{equation}  \label{nabla pk = 1}
\|\nabla p_k \|_{\mathring{M}_{q,\lambda}\left(\Omega\right)}=1
\end{equation}
for all $k$ and
\begin{equation}  \label{sup going to zero}
\sup\left\{ \frac{\int_\Omega \nabla p_k\cdot \nabla \phi}{\|\nabla \phi
\|_{H_{q^{\prime},\lambda}\left(\Omega\right)}}\right\}\to 0.
\end{equation}
By Banach-Alaoglu theorem and by Remark \ref{remark irrotational vector
fields are closed}, we can suppose that $\nabla p_k \rightharpoonup^\ast
\nabla p$ in the weak-star topology of $M_{q,\lambda}\left(\Omega%
\right)=H_{q^{\prime},\lambda}\left(\Omega\right)^\ast$ for some $\nabla p
\in GM_{q,\lambda}\left(\Omega\right)$. Given $\nabla \phi \in
H_{q^{\prime},\lambda}\left(\Omega\right)$, by \eqref{sup going to zero} we
have
\begin{equation*}
\int_{\Omega} \nabla p \cdot \nabla \phi = \lim \int_{\Omega} \nabla p_k
\cdot \nabla \phi = 0.
\end{equation*}
Then, by Lemma \ref{coerciveness lemma}, $p$ is constant in $\Omega$ and
\begin{equation}  \label{nabla pk goes to 0}
\nabla p_k \rightharpoonup^\ast 0
\end{equation}
in the weak-star topology of $M_{q,\lambda}\left(\Omega\right)=H_{q^{%
\prime},\lambda}\left(\Omega\right)^\ast$.

If $\Omega$ is bounded, Theorems \ref{prop Poincare for f L1} and \ref%
{Theorem:RellichKondrachovZorkoSpaces} (Poincar\'{e} inequality and
Rellich-Kondrachov theorem for Zorko spaces) imply that $p_k \in W^1%
\mathring{M}_{q,\lambda}\left(\Omega\right)$ and, up to subsequence, $\| p_k
- \overline{p} \|_{\mathring{M}_{q,\lambda}\left(\Omega\right)}\to 0 $ for
some $\overline{p} \in W^1\mathring{M}_{q,\lambda}\left(\Omega\right)$. By
the weak-star convergence of $\nabla p_k$, we conclude that $\overline{p}$
is constant. Replacing $p_k$ by $p_k - |\Omega|^{-1}\int_{\Omega}p_k$, we
can suppose that $\int_\Omega p_k=0$. Then, $\int_\Omega \overline{p}=\lim
\int_\Omega p_k =0$, which implies that $\overline{p}=0$ and
\begin{equation}  \label{pk goes to 0 bounded}
\|p_k\|_{\mathring{M}_{q,\lambda}\left(\Omega\right)}\to 0.
\end{equation}

Analogously, if $\Omega$ is an exterior domain, we can suppose that $%
\int_{\Omega\cap B}p_k=0$, where $B$ is some sufficiently large open ball
containing $\mathbb{R}^n \backslash \Omega$, and conclude that
\begin{equation}  \label{pk goes to 0 exterior}
\|p_k\|_{\mathring{M}_{q,\lambda}\left(\Omega\cap B\right)}\to 0.
\end{equation}

Now, we split the proof into three parts.

\noindent {\emph{Part 1:}} Let $x_{0}\in \partial \Omega $. After a suitable
rigid movement, we can suppose that $x_{0}=0$ and that $-e_{n}=(0,0,...,-1)$
is the exterior normal to $\Omega $ at $0$. Let $W\subset \mathbb{R}^{n}$ be
a small open neighborhood of $0$ such that there is a function $\sigma \in
C_{0}^{1}(\mathbb{R}^{n-1})$ with $\sigma (0)=0$ and $V:=\Omega \cap W$
coincides with $S_{\sigma }\cap W$, where $S_{\sigma }$ is the bent
half-space $\{x=(x^{\prime },x_{n})\in \mathbb{R}^{n}:x_{n}>\sigma
(x^{\prime })\}$. We observe that $\Vert \nabla \sigma \Vert _{L^{\infty }(%
\mathbb{R}^{n-1})}$ can be supposed arbitrarily small as long $W$ is taken
sufficiently small too. If $\Omega $ is an exterior domain, we can also
suppose that $V\subset \Omega \cap B$.

Let $W^{\prime }\subset W$ be another open neighborhood of $0$ in $\mathbb{R}%
^{n}$ such that $\overline{V^{\prime }}\subset V$ and let $\eta $ be a
smooth cutoff function with $\eta =1$ in $W^{\prime }$ and $\eta =0$ in $%
\mathbb{R}^{n}\backslash W$. Denote $V^{\prime }:=\Omega \cap W^{\prime }$.
We note that
\begin{gather*}
\Vert \nabla (\eta p_{k})\Vert _{\mathring{M}_{q,\lambda }\left( S_{\sigma
}\right) }=\Vert p_{k}\nabla \eta \Vert _{\mathring{M}_{q,\lambda }\left(
V\backslash V^{\prime }\right) }+\Vert \eta \nabla p_{k}\Vert _{\mathring{M}%
_{q,\lambda }\left( V\right) }\leq c\left( \Vert p_{k}\Vert _{\mathring{M}%
_{q,\lambda }\left( V\backslash V^{\prime }\right) }+\Vert \nabla p_{k}\Vert
_{\mathring{M}_{q,\lambda }\left( V\right) }\right) \\
\leq c\left( \Vert p_{k}\Vert _{\mathring{M}_{q,\lambda }\left( V\backslash
V^{\prime }\right) }+1\right) .
\end{gather*}%
In particular, by \eqref{pk goes to 0 bounded} if $\Omega $ is bounded, or %
\eqref{pk goes to 0 exterior} if $\Omega $ is an exterior domain, it follows
that $\{\nabla (\eta p_{k})\}$ is a bounded sequence in $G\mathring{M}%
_{q,\lambda }\left( S_{\sigma }\right) $. Moreover,

\begin{equation}  \label{ineq nabla pk less dk}
\Vert \nabla p_{k}\Vert _{\mathring{M}_{q,\lambda }\left( V^{\prime }\right)
}\leq \Vert \nabla (\eta p_{k})\Vert _{\mathring{M}_{q,\lambda }\left(
S_{\sigma }\right) }\leq c\sup \left\{ \frac{\int_{S_{\sigma }}\nabla (\eta
p_{k})\cdot \nabla \phi }{\Vert \nabla \phi \Vert _{H_{q^{\prime },\lambda
}\left( S_{\sigma }\right) }}\right\} ,
\end{equation}%
where the supremum is taken over all non null functions $\nabla \phi \in
GH_{q^{\prime },\lambda }\left( S_{\sigma }\right) $. For each $k$, let us
denote
\begin{equation*}
d_{k}:=\sup \left\{ \frac{\int_{S_{\sigma }}\nabla (\eta p_{k})\cdot \nabla
\phi }{\Vert \nabla \phi \Vert _{H_{q^{\prime },\lambda }\left( S_{\sigma
}\right) }}\right\}
\end{equation*}%
and let $\nabla \phi _{k}\in GH_{q^{\prime },\lambda }\left( S_{\sigma
}\right) $ such that $\Vert \nabla \phi _{k}\Vert _{H_{q^{\prime },\lambda
}\left( S_{\sigma }\right) }=1$ and
\begin{equation*}
d_{k}\leq 1/k+\int_{S_{\sigma }}\nabla (\eta p_{k})\cdot \nabla \phi
_{k}=1/k+\int_{V}\nabla (\eta p_{k})\cdot \nabla \phi _{k}.
\end{equation*}%
We can suppose that $\int_{V\backslash V^{\prime }}\phi _{k}=0$. Then, by
Theorem \ref{theorem Poincare inequality block} (Poincar\'{e} inequality for
block spaces), item (III), we have that $\phi _{k}\in W^{1}H_{q^{\prime
},\lambda }\left( V\backslash V^{\prime }\right) $ and
\begin{equation*}
\Vert \phi _{k}\Vert _{H_{q^{\prime },\lambda }\left( V\backslash V^{\prime
}\right) }\leq c\Vert \nabla \phi _{k}\Vert _{H_{q^{\prime },\lambda }\left(
V\backslash V^{\prime }\right) }\leq c.
\end{equation*}%
Therefore, $\{\phi _{k}\}$ is a bounded sequence in $W^{1}H_{q^{\prime
},\lambda }\left( V\backslash V^{\prime }\right) $ and by Theorem \ref%
{compact embedding block} (Rellich-Kondrachov for block spaces) there is $%
\phi \in W^{1}H_{q^{\prime },\lambda }\left( V\backslash V^{\prime }\right) $
such that, up to subsequence,
\begin{equation}
\Vert \phi _{k}-\phi \Vert _{H_{q^{\prime },\lambda }\left( V\backslash
V^{\prime }\right) }\rightarrow 0  \label{phik converges}
\end{equation}%
and
\begin{equation}
\nabla \phi _{k}\rightharpoonup ^{\ast }\nabla \phi
\label{nabla phik converges}
\end{equation}%
in the weak-star topology of $H_{q^{\prime },\lambda }\left( V\backslash
V^{\prime }\right) =\mathring{M}_{q,\lambda }\left( V\backslash V^{\prime
}\right) ^{\ast }$. Then,
\begin{gather*}
d_{k}-1/k\leq \int_{V\backslash V^{\prime }}\nabla (\eta p_{k})\cdot \nabla
\phi _{k}=\int_{V\backslash V^{\prime }}p_{k}\nabla \eta \cdot \nabla \phi
_{k}+\int_{V}\eta \nabla p_{k}\cdot \nabla \phi _{k} \\
=\int_{V\backslash V^{\prime }}p_{k}\nabla \eta \cdot \nabla \phi
_{k}+\int_{V}\nabla p_{k}\cdot \nabla (\eta \phi _{k})-\int_{V\backslash
V^{\prime }}\phi _{k}\nabla p_{k}\cdot \nabla \eta .
\end{gather*}%
All the three integrals above converge to zero: the first one, by
\eqref{pk
goes to 0 bounded} (or \eqref{pk goes to 0 exterior}) and
\eqref{nabla phik
converges}; the second one, by \eqref{sup going to zero} and since $\{\nabla
(\eta \phi _{k})\}$ is bounded in $H_{q^{\prime },\lambda }\left( \Omega
\right) $; and the third one, by \eqref{nabla pk goes to 0} and
\eqref{phik
converges}. Therefore, $d_{k}\rightarrow 0$ and, by
\eqref{ineq nabla pk
less dk},
\begin{equation*}
\Vert \nabla p_{k}\Vert _{\mathring{M}_{q,\lambda }\left( V^{\prime }\right)
}\rightarrow 0.
\end{equation*}

By compactness, we can cover $\partial \Omega$ with a finite number of open
sets like $W^{\prime}$ and then conclude that there is a open set $%
W_0\supset \partial \Omega$ (formed by the finite union of such sets $%
W^{\prime}$) such that for $V_0:=\Omega\cap W_0$, we have
\begin{equation}  \label{estimate boundary}
\|\nabla p_k\|_{\mathring{M}_{q,\lambda}\left(V_0\right)} \to 0.
\end{equation}
If $\Omega$ is an exterior domain, we also can suppose that $V_0 \subset
\Omega\cap B$.

\noindent{\emph{Part 2:}} The part 2 consists in to prove an estimate such
as \eqref{estimate boundary} for the interior of $\Omega$. Let $U$ be an
open set such that $\overline{U}\subset \Omega$ and $\overline{U^{\prime}}%
\subset U$, where $U^{\prime}:=\Omega \backslash \overline{V_0}$. The
argument is very similar to that in the part 1 with $V, V^{\prime}$ replaced
by $U,U^{\prime}$.

Let $\theta$ be a smooth cutoff function with $\theta =1$ in $U^{\prime}$
and $\theta=0$ in $\mathbb{R}^n \backslash U$. Then, $\{\nabla (\theta
p_k)\} $ is a bounded sequence in $G\mathring{M}_{q,\lambda}\left(\mathbb{R}%
^n\right)$. Moreover,
\begin{align}  \label{ineq nabla pk less dk 2}
\|\nabla p_k\|_{\mathring{M}_{q,\lambda}\left(U^{\prime}\right)} \leq
\|\nabla (\theta p_k)\|_{\mathring{M}_{q,\lambda}\left(\mathbb{R}^n\right)}
\leq c \sup\left\{ \frac{\int_{\mathbb{R}^n} \nabla (\theta p_k)\cdot \nabla
\phi}{\|\nabla \phi \|_{H_{q^{\prime},\lambda}\left(\mathbb{R}^n\right)}}%
\right\},
\end{align}
where the supremum is taken over all non null functions $\nabla\phi \in
GH_{q^{\prime},\lambda}\left(\mathbb{R}^n\right)$. For each $k$, let us
denote
\begin{equation*}
b_k:=\sup\left\{ \frac{\int_{\mathbb{R}^n} \nabla (\theta p_k)\cdot \nabla
\phi}{\|\nabla \phi \|_{H_{q^{\prime},\lambda}\left(\mathbb{R}^n\right)}}%
\right\}
\end{equation*}
and let $\nabla\phi_k \in GH_{q^{\prime},\lambda}\left(\mathbb{R}^n\right)$
such that $\|\nabla \phi_k\|_{H_{q^{\prime},\lambda}\left(\mathbb{R}%
^n\right)}=1$ and
\begin{equation*}
b_k \leq 1/k + \int_{\mathbb{R}^n}\nabla(\theta p_k)\cdot \nabla \phi_k =
1/k + \int_{U}\nabla(\theta p_k)\cdot \nabla \phi_k.
\end{equation*}
We can suppose that $\int_{U\backslash U^{\prime}}\phi_k=0$. Then, $%
\{\phi_k\}$ is a bounded sequence in $W^1H_{q^{\prime},\lambda}\left(U%
\backslash U^{\prime}\right)$ and there is $\phi \in
W^1H_{q^{\prime},\lambda}\left(U\backslash U^{\prime}\right)$ such that, up
to subsequence,
\begin{equation}  \label{phik converges 2}
\|\phi_k- \phi\|_{H_{q^{\prime},\lambda}\left(U\backslash
U^{\prime}\right)}\to 0
\end{equation}
and
\begin{equation}  \label{nabla phik converges 2}
\nabla \phi_k \rightharpoonup^\ast \nabla \phi
\end{equation}
in the weak-star topology of $H_{q^{\prime},\lambda}\left(U\backslash
U^{\prime}\right)=\mathring{M}_{q,\lambda}\left(U\backslash
U^{\prime}\right)^\ast$. Then,
\begin{align*}
b_k-1/k & \leq \int_{U\backslash U^{\prime}} p_k\nabla \theta\cdot \nabla
\phi_k + \int_{U} \nabla p_k\cdot \nabla (\theta\phi_k) - \int_{U\backslash
U^{\prime}} \phi_k\nabla p _k \cdot \nabla \theta.
\end{align*}
As in the part 1, the three integrals above converge to zero. Therefore, $%
b_k \to 0$ and, by \eqref{ineq nabla pk less dk 2},
\begin{equation}  \label{estimate interior}
\|\nabla p_k\|_{\mathring{M}_{q,\lambda}\left(U^{\prime}\right)} \to 0.
\end{equation}

\noindent{\emph{Part 3:}} By \eqref{estimate boundary} and
\eqref{estimate
interior} we conclude that
\begin{equation}  \label{global estimate}
\|\nabla p_k\|_{\mathring{M}_{q,\lambda}\left(\Omega\right)} \to 0,
\end{equation}
which contradicts \eqref{nabla pk = 1}. \fin

As a consequence of Theorem \ref{variational inequality theorem}, we obtain
the Helmholtz decomposition for $\mathring{M}_{q,\lambda }\left( \Omega
\right) $ with $\Omega $ being either a bounded or an exterior $C^{1}$
domain.

\begin{theorem}
\label{Theorem:HelmholtzDecompositionZorkoSpacesBoundedExteriorDomain} Let $%
\Omega$ be a bounded or exterior $C^1$ domain. Then the Helmholtz
decomposition of $\mathring{M}_{q,\lambda}\left(\Omega\right)$ holds.
\end{theorem}

\noindent{\textbf{Proof.}} Let us consider the bounded linear operator
\begin{equation*}
T:G\mathring{M}_{q,\lambda}\left(\Omega\right)\to
GH_{q^{\prime},\lambda}\left(\Omega\right)^\ast
\end{equation*}
given by
\begin{equation*}
T(\nabla f)(\nabla g)=\int_{\Omega}\nabla f \cdot \nabla g,
\end{equation*}
for $\nabla f \in G\mathring{M}_{q,\lambda}\left(\Omega\right), \nabla g \in
GH_{q^{\prime},\lambda}\left(\Omega\right).$ Then, Theorem \ref{variational
inequality theorem} can be restated as
\begin{equation}  \label{coerciveness restated}
\|\nabla p\|_{\mathring{M}_{q,\lambda}\left(\Omega\right)}\leq c \|T(\nabla
p)\|_{GH_{q^{\prime},\lambda}\left(\Omega\right)^\ast},
\end{equation}
which implies, in particular, that $T$ has closed range. Moreover, by the
Helmholtz decomposition in Lebesgue spaces, the restriction $%
S:=T|_{GL^{n/\alpha}(\Omega)}$, where $\alpha:=(n-\lambda)/q$, is an
isomorphism between $GL^{n/\alpha}(\Omega)$ and $GL^{n/(n-\alpha)}(\Omega)^%
\ast$ (see Theorem 6.1 in \cite{simader1992new} or Theorem 4 in \cite%
{frohlich2000helmholtz}). In particular, the adjoint operator $S^\ast$ of $S$
is an isomorphism between $GL^{n/(n-\alpha)}(\Omega)^{\ast\ast}$ to $%
GL^{n/\alpha}(\Omega)^\ast$ and, since $T^\ast$ is the restriction of $%
S^\ast $ to $GH_{q^{\prime},\lambda}\left(\Omega\right)^{\ast\ast}$, we have
that $T^\ast$ is injective. Then, by the closed range theorem, $T$ is
surjective, which along with \eqref{coerciveness restated}, means that $T$
is an isomorphism.

Finally, given $\mathbf{u} \in \mathring{M}_{q,\lambda}\left(\Omega\right)$,
consider the functional $\mathcal{F} \in
GH_{q^{\prime},\lambda}\left(\Omega\right)^\ast$ defined by $\nabla \phi
\mapsto \int_{\Omega}\mathbf{u}\cdot \nabla \phi$. As we have proven, there
is a unique $\nabla p \in G\mathring{M}_{q,\lambda}\left(\Omega\right)$ such
that $T(\nabla p)=\mathcal{F}$, that is,
\begin{equation*}
\int_{\Omega} \nabla p \cdot \nabla \phi = \int_\Omega \mathbf{u}\cdot
\nabla \phi,
\end{equation*}
for all $\nabla \phi \in GH_{q^{\prime},\lambda}\left(\Omega\right)$. Then,
by Lemma \ref{lemma de equivalencia}, the Helmholtz decomposition for $%
\mathring{M}_{q,\lambda}\left(\Omega\right)$ holds. Observe that
\begin{equation*}
\|\nabla p \|_{\mathring{M}_{q,\lambda}\left(\Omega\right)}\leq c \|T(\nabla
p)\|_{GH_{q^{\prime},\lambda}\left(\Omega\right)^\ast}=c\|\mathcal{F}%
\|_{GH_{q^{\prime},\lambda}\left(\Omega\right)^\ast} \leq c \|\mathbf{u}\|_{%
\mathring{M}_{q,\lambda}\left(\Omega\right)}.
\end{equation*}
\fin

\subsection{Helmholtz Decomposition for $M_{q,\protect\lambda%
}\left(\Omega\right)$ and $H_{q,\protect\lambda}(\Omega)$}

\label{subsection Helmholtz Decomposition for Morrey and Block spaces}

Finally, we consider the Helmholtz decomposition for $M_{q,\lambda}\left(%
\Omega\right)$ and $H_{q,\lambda}(\Omega)$, actuality, for the sake of
convenience, for $M_{q,\lambda}\left(\Omega\right)$ and $H_{q^{\prime},%
\lambda}\left(\Omega\right)$. We remind that $\Omega$ can be a bounded or
exterior domain with $C^1$ boundary, $\mathbb{R}^n$ or $\mathbb{R}^n_+$. We
argument by duality. For this purpose, it is useful to consider the
Helmholtz projection operator $\mathbf{P}_{\mathring{M}_{q,\lambda}\left(%
\Omega\right)}:\mathring{M}_{q,\lambda}\left(\Omega\right)\to\mathring{M}%
_{q,\lambda}\left(\Omega\right)$, which we abbreviate $\mathbf{P}=\mathbf{P}%
_{\mathring{M}_{q,\lambda}\left(\Omega\right)}$, with range $R(\mathbf{P})=S%
\mathring{M}_{q,\lambda}\left(\Omega\right)$ and kernel $N(\mathbf{P})=G%
\mathring{M}_{q,\lambda}\left(\Omega\right)$. Then $\mathbf{P^{\ast}}$ and $%
\mathbf{P^{\ast\ast}}$ are bounded projection operators in $%
H_{q^{\prime},\lambda}\left(\Omega\right)$ and $M_{q,\lambda}\left(\Omega%
\right)$, so
\begin{equation*}
H_{q^{\prime},\lambda}\left(\Omega\right)=N(\mathbf{P^{\ast}})\oplus R(%
\mathbf{P^{\ast}})\text{ and }M_{q,\lambda}\left(\Omega\right)=N(\mathbf{%
P^{\ast\ast}})\oplus R(\mathbf{P^{\ast\ast}}).
\end{equation*}

By properties of annihilators,
\begin{equation*}
N(\mathbf{P^{\ast}})=R(\mathbf{P})^\perp=S\mathring{M}_{q,\lambda}\left(%
\Omega\right)^\perp.
\end{equation*}
From Lemma \ref{lema de Rham} (de Rham), we obtain that $S\mathring{M}%
_{q,\lambda}\left(\Omega\right)^\perp=GH_{q^{\prime},\lambda}\left(\Omega%
\right)$. Then,
\begin{equation*}
N(\mathbf{P^{\ast}})=GH_{q^{\prime},\lambda}\left(\Omega\right).
\end{equation*}
On the other hand,
\begin{equation*}
R(\mathbf{P^{\ast}})=N(\mathbf{P})^\perp = G\mathring{M}_{q,\lambda}\left(%
\Omega\right)^\perp.
\end{equation*}
If we consider the space $C^\infty_{0,\sigma}(\Omega) \subset
H_{q^{\prime},\lambda}\left(\Omega\right)$ as a subset of $\mathring{M}%
_{q,\lambda}\left(\Omega\right)^\ast$, then, also from Lemma \ref{lema de
Rham}, its preannihilator $^\perp C^\infty_{0,\sigma}(\Omega) \subset%
\mathring{M}_{q,\lambda}\left(\Omega\right)$ coincides with $G\mathring{M}%
_{q,\lambda}\left(\Omega\right)$. Therefore, $R(\mathbf{P^{\ast}})= \left(
^\perp C^\infty_{0,\sigma}(\Omega) \right)^\perp=$ weak-star closure of $%
C^\infty_{0,\sigma}(\Omega)$ in $\mathring{M}_{q,\lambda}\left(\Omega%
\right)^\ast$, that is,
\begin{align*}
R(\mathbf{P^{\ast}})=SH_{q^{\prime},\lambda}\left(\Omega\right).
\end{align*}
Then the decomposition holds for $H_{q^{\prime},\lambda}\left(\Omega\right)$%
:
\begin{equation*}
H_{q^{\prime},\lambda}\left(\Omega\right)=GH_{q^{\prime},\lambda}\left(%
\Omega\right)\oplus SH_{q^{\prime},\lambda}\left(\Omega\right).
\end{equation*}

Similarly, for $M_{q,\lambda}\left(\Omega\right)$,
\begin{align*}
R(\mathbf{P^{\ast\ast}})=N(\mathbf{P^{\ast}})^\perp
=GH_{q^{\prime},\lambda}\left(\Omega\right)^\perp = \text{ weak* closure of }%
C^\infty_{0,\sigma}(\Omega) \text{ in } H_{q^{\prime},\lambda}\left(\Omega%
\right)^\ast = SM_{q,\lambda}\left(\Omega\right).
\end{align*}
Also,
\begin{equation*}
N(\mathbf{P^{\ast\ast}})=R(\mathbf{P^{\ast}})^\perp =
SH_{q^{\prime},\lambda}\left(\Omega\right)^\perp \subset
GM_{q,\lambda}\left(\Omega\right).
\end{equation*}
Therefore,
\begin{equation*}
N(\mathbf{P^{\ast\ast}})\oplus R(\mathbf{P^{\ast\ast}})=GM_{q,\lambda}\left(%
\Omega\right)+ SM_{q,\lambda}\left(\Omega\right).
\end{equation*}
Since $SM_{q,\lambda}\left(\Omega\right)\cap
GM_{q,\lambda}\left(\Omega\right)\subset S L_w^q(\Omega) \cap G
L_w^q(\Omega)=\{\mathbf{0}\}$, we conclude that $N(\mathbf{P^{\ast\ast}}%
)=GM_{q,\lambda}\left(\Omega\right)$ and the Helmholtz decomposition
\begin{equation*}
M_{q,\lambda}\left(\Omega\right)=GM_{q,\lambda}\left(\Omega\right)\oplus
SM_{q,\lambda}\left(\Omega\right)
\end{equation*}
holds.

\begin{remark}
\label{remark weak-star closure} The spaces $SM_{q,\lambda }(\Omega )$ and $%
SH_{q,\lambda }(\Omega )$ are defined as the closures of $C_{0,\sigma
}^{\infty }(\Omega )$ in the \emph{weak-star} topology of $H_{q,\lambda
}(\Omega )^{\ast }$ and $\mathring{M}_{q,\lambda }(\Omega )^{\ast }$,
respectively, and not in the norm topology, as is usual for other classes of
function spaces. This is an effect of the non-reflexivity of Morrey (Zorko
or block) spaces. We observe that, for a reflexive function space $X$--such
as a weighted Lebesgue space, for instance--the weak-star closure of $%
C_{0,\sigma }^{\infty }\cap X$ coincides with the weak closure, which, since
$C_{0,\sigma }^{\infty }$ is convex, also coincides with the strong closure.
Moreover, even for $L^{p}(\Omega )$, a more relaxed definition of $%
SL^{q}(\Omega )$ rather than the usual one can be necessary; for instance,
if one considers more general types of domains (see \cite%
{simader2014necessary} and references therein). We note that, as expected, $%
SM_{q,\lambda }\left( \Omega \right) $ is a strictly larger space than $S%
\mathring{M}_{q,\lambda }\left( \Omega \right) $. To see this, consider, for
instance, the vector field $\mathbf{v(x)=}|x-x^{\ast }|^{-1-\alpha }\mathbf{u%
}(x-x^{\ast }),$ $x\in \Omega ,$ where $x^{\ast }\in \Omega ,$ $\alpha
=(n-\lambda )/q$, and $\mathbf{u}(x)=(x_{2},-x_{1},0,\dots ,0)$. Note that $%
\mathbf{v\in }SM_{q,\lambda }\left( \Omega \right) $ but $\mathbf{v}\notin S%
\mathring{M}_{q,\lambda }\left( \Omega \right) $. However, $GM_{q,\lambda
}\left( \Omega \right) =G\mathring{M}_{q,\lambda }\left( \Omega \right) $.
Indeed, as a consequence of the proof of Theorem \ref%
{Theorem:HelmholtzDecompositionZorkoSpacesBoundedExteriorDomain}, we have $G%
\mathring{M}_{q,\lambda }\left( \Omega \right) =GH_{q^{\prime },\lambda
}\left( \Omega \right) ^{\ast }$ (in the sense of duality induced by
integration). On the other hand, as seen above, $GM_{q,\lambda }\left(
\Omega \right) =SH_{q^{\prime },\lambda }\left( \Omega \right) ^{\perp }$.
Since $GH_{q^{\prime },\lambda }\left( \Omega \right) ^{\ast }=SH_{q^{\prime
},\lambda }\left( \Omega \right) ^{\perp }$, we have
\begin{equation*}
GM_{q,\lambda }\left( \Omega \right) =G\mathring{M}_{q,\lambda }\left(
\Omega \right) .
\end{equation*}
\end{remark}

\

\noindent \textbf{{\large {Acknowledgments.}}} L. C. F. Ferreira was
supported by CNPq (grant 308799/2019-4), Brazil. M. G. Santana was supported
by CNPq (grant 140395/2021-0), Brazil.

\end{document}